%% file: unrooted_maps.tex
\documentclass{amsart}

\usepackage{latexsym, amssymb, amsmath}
\usepackage{epsfig}
\usepackage{graphics}
\usepackage{color}
\usepackage{subfigure,pstricks,pst-node}

\input{Figures/head.tex}

\newtheorem{theorem}{Theorem}
\newtheorem{lemma}[theorem]{Lemma}

\newtheorem{proposition}[theorem]{Proposition}

\newcommand{\tb}{w_{\bullet}}
\newcommand{\tw}{w_{\circ}}

\newcommand{\xb}{x_{\bullet}}
\newcommand{\xw}{x_{\circ}}
\newcommand{\yb}{y_{\bullet}}
\newcommand{\yw}{y_{\circ}}
\newcommand{\zb}{z_{\bullet}}
\newcommand{\zw}{z_{\circ}}
\newcommand{\ovw}{W(y)/y}
\newcommand{\ovg}{J(y)/y}
\newcommand{\Gvv}{\mathcal{G}_{vv}^{(2)}}
\newcommand{\gvv}{G_{vv}^{(2)}}
\newcommand{\Qvv}{\mathcal{H}_{vv}^{(2)}}
\newcommand{\qvv}{H_{vv}^{(2)}}
\newcommand{\Gfvtt}{\mathcal{G}_{v\mathfrak{f}'}^{(2)}}
\newcommand{\gfvtt}{G_{v\mathfrak{f}'}^{(2)}}
\newcommand{\Qfvtt}{\mathcal{H}_{v\mathfrak{f}'}^{(2)}}
\newcommand{\qfvtt}{H_{v\mathfrak{f}'}^{(2)}}
\newcommand{\Gfvt}{\mathcal{G}_{v\mathfrak{f}}^{(2)}}
\newcommand{\gfvt}{G_{v\mathfrak{f}}^{(2)}}
\newcommand{\Qfvt}{\mathcal{H}_{v\mathfrak{f}}^{(2)}}
\newcommand{\qfvt}{H_{v\mathfrak{f}}^{(2)}}

\newcommand{\Gfftt}{\mathcal{G}_{\mathfrak{f}\mathfrak{f}'}^{(2)}}
\newcommand{\gfftt}{G_{\mathfrak{f}\mathfrak{f}'}^{(2)}}
\newcommand{\Qfftt}{\mathcal{H}_{\mathfrak{f}\mathfrak{f}'}^{(2)}}
\newcommand{\qfftt}{H_{\mathfrak{f}\mathfrak{f}'}^{(2)}}
\newcommand{\Gfft}{\mathcal{G}_{\mathfrak{f}\mathfrak{f}}^{(2)}}
\newcommand{\gfft}{G_{\mathfrak{f}\mathfrak{f}}^{(2)}}
\newcommand{\Qfft}{\mathcal{H}_{\mathfrak{f}\mathfrak{f}}^{(2)}}
\newcommand{\qfft}{H_{\mathfrak{f}\mathfrak{f}}^{(2)}}
\newcommand{\gfbtt}{G_{b\mathfrak{f}'}^{(2)}}
\newcommand{\gfwtt}{G_{w\mathfrak{f}'}^{(2)}}
\newcommand{\gfbt}{G_{b\mathfrak{f}}^{(2)}}
\newcommand{\gfwt}{G_{w\mathfrak{f}}^{(2)}}
\newcommand{\gbb}{G_{bb}^{(2)}}
\newcommand{\gbw}{G_{bw}^{(2)}}
\newcommand{\gww}{G_{ww}^{(2)}}
\newcommand{\yW}{\left( W/\yw,W/\yb\right)}
\newcommand{\qfbt}{H_{b\mathfrak{f}}^{(2)}}
\newcommand{\Qfbt}{\mathcal{H}_{b\mathfrak{f}}^{(2)}}
\newcommand{\Gfbt}{\mathcal{G}_{b\mathfrak{f}}^{(2)}}
\newcommand{\Qfwt}{\mathcal{H}_{w\mathfrak{f}}^{(2)}}
\newcommand{\Gfwt}{\mathcal{G}_{w\mathfrak{f}}^{(2)}}
\newcommand{\qfwt}{H_{w\mathfrak{f}}^{(2)}}
\newcommand{\qfbtt}{H_{b\mathfrak{f}'}^{(2)}}
\newcommand{\Qfbtt}{\mathcal{H}_{b\mathfrak{f}'}^{(2)}}
\newcommand{\Gfbtt}{\mathcal{G}_{b\mathfrak{f}'}^{(2)}}
\newcommand{\Gfwtt}{\mathcal{G}_{w\mathfrak{f}'}^{(2)}}
\newcommand{\Qfwtt}{\mathcal{H}_{w\mathfrak{f}'}^{(2)}}
\newcommand{\qfwtt}{H_{w\mathfrak{f}'}^{(2)}}
\newcommand{\qbb}{H_{bb}^{(2)}}
\newcommand{\Qbb}{\mathcal{H}_{bb}^{(2)}}
\newcommand{\Gbb}{\mathcal{G}_{bb}^{(2)}}
\newcommand{\qww}{H_{ww}^{(2)}}
\newcommand{\qbw}{H_{bw}^{(2)}}
\newcommand{\Qbw}{\mathcal{H}_{bw}^{(2)}}
\newcommand{\Gbw}{\mathcal{G}_{bw}^{(2)}}
\newcommand{\Qww}{\mathcal{H}_{ww}^{(2)}}
\newcommand{\Gww}{\mathcal{G}_{ww}^{(2)}}

\title{Counting unrooted maps using tree-decomposition}

\author[\'E. Fusy]{\'Eric Fusy}
\address{\'Eric Fusy, INRIA Rocquencourt, Projet ALGO
BP 105, 78153 Le Chesnay, and \'Ecole Polytechnique, LIX}
\email{Eric.Fusy@inria.fr}

\begin{document}
\bibliographystyle{alpha}

\begin{abstract}
We present a new method to count unrooted maps on the sphere up to
orientation-preserving homeomorphisms. The principle, called 
\emph{tree-decomposition}, is to deform a map into an arborescent 
structure whose nodes are occupied by constrained maps. 
Tree-decomposition turns out to be very efficient and flexible for 
the enumeration of constrained families of maps. In this article, 
the method is applied to count unrooted 2-connected maps and, more 
importantly, to count unrooted 3-connected maps, which correspond 
to the combinatorial types of oriented convex polyhedra. Our
method improves significantly on the previously best-known complexity 
to enumerate
unrooted 3-connected maps.
\end{abstract}

\maketitle

\noindent\textbf{Acknowledgments.} The author would like to thank
Gilles Schaeffer for his invaluable help in developing this new
method. In particular he pointed out the idea of tree-decomposition and
helped to do some calculations and to correct the article.

\section*{Introduction}
The enumeration of unrooted maps has been a well-studied problem for
more than 20 years. The first general method for the enumeration of unrooted maps on the sphere up
to orientation-preserving homeomorphisms was developed by Liskovets~\cite{Li81}. It is based on two main
tools: Burnside's formula and study of the structure of the quotient maps.

With an adaptation of Burnside's (orbit counting) formula, counting unrooted maps comes down to counting rooted maps with
a symmetry of rotation. For a family of maps enumerated
according to the number $n$ of edges, we write
respectively $c_n$, $c_n'$ and $c_n^{(k)}$ for the number of unrooted
maps, rooted maps and rooted maps with a symmetry of order $k\geq 2$; then
$c_n$ can be computed with the formula:

\begin{equation}
\label{eq:burnsi}
c_n=\frac{1}{2n}\left(c_n'+\sum_{k=2}^n\phi(k)c_n^{(k)}\right)
\end{equation}
and a similar formula exists for the enumeration according to the
number of vertices and faces, see Section~\ref{sec:enum}. We represent rooted 
maps with a symmetry of order $k\geq 2$ as
$k$-rooted maps, which are maps with $k$ undistinguishable
roots. The \emph{quotient map} of such a symmetric map is
 a rooted map with two marked cells, the cells intersected by the rotation-axis, that are either a vertex or the
middle of a face or the middle of an edge. The enumeration of these maps is 
easy
to handle for the family of unconstrained maps~\cite{Li81}. 
The method of quotient-maps can also be adapted 
to the enumeration of some
families of constrained maps, such as loopless maps
~\cite{LW04}, eulerian and unicursal maps~\cite{LW02} and 2-connected
maps~\cite{LW83} but the the structure of the quotient maps is less easy to characterize and to handle for these families.

\begin{figure}
\begin{center}
\includegraphics[width=14cm]{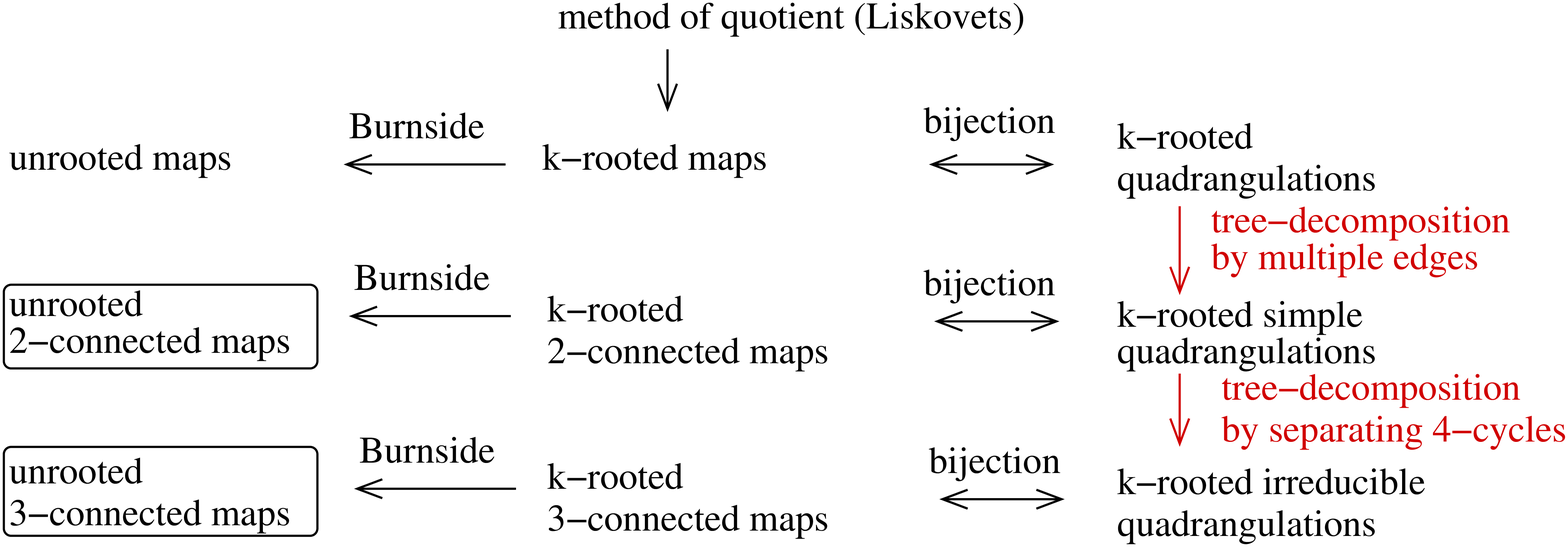}
\end{center}
\caption{The scheme of the method to enumerate unrooted 2-connected
and unrooted 3-connected maps}
\label{figure:diagschema}
\end{figure}

\begin{figure}
  \begin{center}
    \includegraphics[width=13cm]{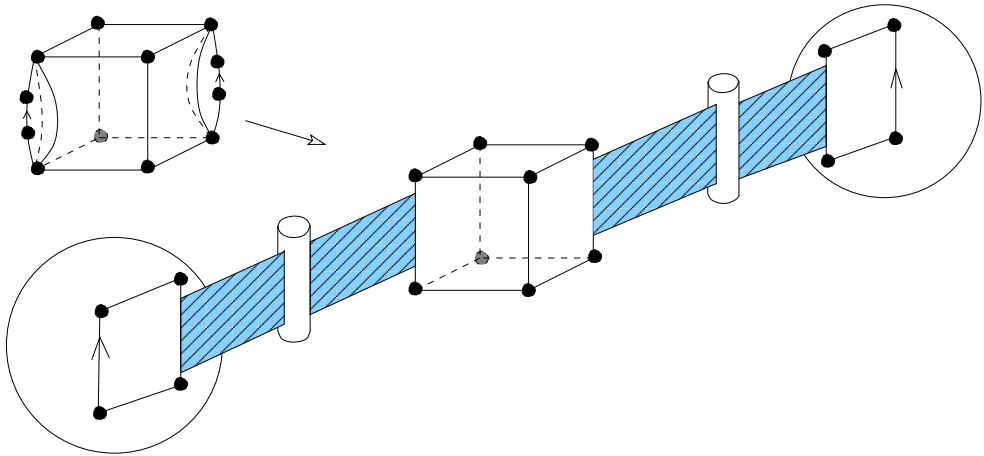}
  \end{center}
  \caption{Repercussion of the symmetry of a $k$-rooted
quadrangulation on its decomposition-tree.}
\label{figure:repsymqu}
\end{figure}

In this article, we introduce a new general method for the enumeration of 
unrooted 
maps of a constrained family, based on the concept of tree-decomposition. 
We apply
the method to the enumeration of unrooted 2-connected
and, above all, of unrooted 3-connected maps, already counted by
Walsh~\cite{Wa82}, but with a costly step of extraction of coefficients. 
In order to apply the method of tree-decomposition for 2-connected and 
3-connected maps, we prefer to work with quadrangulations rather than with 
maps. Indeed, a well-known bijection between maps and quadrangulations, 
recalled in Section~\ref{sec:bijkrooted}, ensures that counting 2-connected maps and 
3-connected maps is respectively equivalent to counting \emph{simple} 
quadrangulations (i.e. quadrangulations without multiple edges) and 
\emph{irreducible} quadrangulations (i.e. quadrangulations without 
separating 4-cycles). Then, we 
introduce two tree-decompositions on quadrangulations.   A first 
tree-decomposition ``by multiple edges'', ensures that a quadrangulation
can be seen as an arborescent structure with nodes that are simple 
quadrangulations. The symmetry of order $k$ of a $k$-rooted quadrangulation 
fixes the decomposition-tree, hence it also fixes the node at the centre of 
the tree, called the \emph{core-node}, see Figure~\ref{figure:repsymqu} that best summarizes 
the essence of the method. This yields an equation linking the
generating function of $k$-rooted simple quadrangulations and the generating
functions  of $k$-rooted quadrangulations, which are easy to obtain using the
method of quotient map (see Section~\ref{sec:algunconstr} where we 
briefly re-derive the results of Liskovets~\cite{Li81} for generating 
functions, and Sections~\ref{sec:enum_edges_1conn} 
and~\ref{sec:enum_vertices_1conn} for all explicit expressions). 
Then, a second tree-decomposition ``by separating 
4-cycles'', introduced by Kunz-Jacques and Schaeffer in~\cite{Ku} 
for the enumeration of prime alternating links, states that a 
simple quadrangulation can be seen as an arborescent 
structure with nodes that are irreducible quadrangulations. 
In a similar way as 
for the first tree-decomposition, the symmetry of a $k$-rooted simple 
quadrangulation also fixes the ``core'' of the decomposition-tree, yielding 
equations linking the generating
functions of $k$-rooted irreducible quadrangulations and the generating 
functions of
$k$-rooted simple quadrangulations, which have already been obtained thanks 
to the
first tree-decomposition. However, for this second tree-decomposition, 
a careful treatment of cases has to be done 
in Section~\ref{sec:2rootedirreducible} for 2-rooted objects, 
since the core of the decomposition-tree can be different from the 
centre of the tree. Once the generating functions of $k$-rooted simple 
quadrangulations (equal to those of $k$-rooted 2-connected maps) and the 
generating functions of $k$-rooted irreducible quadrangulations (equal to 
those of $k$-rooted 3-connected maps) are obtained, Burnside's formula 
(\ref{eq:burnsi}) yields respectively the enumeration of unrooted 2-connected 
and unrooted 3-connected maps, see Figure~\ref{figure:diagschema} for a 
summarizing diagram.

\textbf{Main results.}
Two results are obtained: a theorem about the
\emph{algebraic structure} of $k$-rooted maps and a theorem giving the
complexity of enumeration of unrooted 2-connected and unrooted
3-connected maps. First we need a few notations. Given a series $\alpha(t)$, 
a series $f(t)$ is said
$\alpha$-rational if there exists a rational function $R(T)$ such that
$f(t)=R(\alpha(t))$. Given two series in two variables
$\alpha_1(t_{\bullet},t_{\circ})$ and
$\alpha_2(t_{\bullet},t_{\circ})$, a series in two variables 
$f(t_{\bullet},t_{\circ})$ is said $(\alpha_1,\alpha_2)$-rational if there 
exists a
rational function $R(T_1,T_2)$ in two variables such that
$f(t_{\bullet},t_{\circ})=R(\alpha_1(t_{\bullet},t_{\circ}),\alpha_2(t_{\bullet},t_{\circ}))$.

Now we introduce the three ``easily'' algebraic series in one variable
(they correspond to families of trees) $\beta(x)$, $\eta(y)$ and
$\gamma(z)$~\footnote{We use three different variable names $x$, $y$,
$z$, because they will later be linked by relations of change of variable.} 
given by

$$
\beta(x)=x+3\beta(x)^2 ,\ \ \ \ \ \ \ \ \ \ \ \ \ \ 
\eta(y)=\frac{y}{(1-\eta(y))^2} ,\ \ \ \ \ \ \ \ \ \ \ \ \ 
\gamma(z)=z(1+\gamma(z))^2
$$
and their versions in two variables $\beta_{1,2}(\xb,\xw)$,
$\eta_{1,2}(\yb,\yw)$, and $\gamma_{1,2}(\zb,\zw)$ (corresponding to
bicolored trees of the respective families) given by

\begin{minipage}{3.8cm}
$$
\left\{
\begin{array}{rcl}
\beta_1&=&\xb +\beta_1^2+2\beta_1\beta_2 \\
\beta_2&=&\xw +\beta_2^2+2\beta_1\beta_2
\end{array}
\right. ,
$$
\end{minipage}
\begin{minipage}{5.2cm}
$$
\left\{
\begin{array}{rcl}
\eta_1&=&\yb/(1-\eta_2)^2\\
\eta_2&=&\yw/(1-\eta_1)^2
\end{array}
\right. ,
$$
\end{minipage}
 \begin{minipage}{3.7cm}
$$
\left\{
\begin{array}{rcl}
\gamma_1&=&\zb(1+\gamma_2)^2\\
\gamma_2&=&\zw(1+\gamma_1)^2
\end{array}
\right. .
$$
\end{minipage}

It is already known that the generating function of rooted maps is 
$\beta$-rational in one variable~\cite{Tu63} and $(\beta_1,\beta_2)$-rational 
in two variables~\cite{Sc97}; that the generating function of rooted 2-connected maps 
is $\eta$-rational in one variable and $(\eta_1,\eta_2)$-rational in two 
variables~\cite{BT}; and  that the generating function of rooted 3-connected planar 
maps is $\gamma$-rational in one variable and $(\gamma_1,\gamma_2)$-rational 
in two variables~\cite{MS68}. The following theorem states that the same property also 
holds for $k$-rooted maps:
\begin{theorem}
\label{theo:first}
For $k\geq 2$, the series of $k$-rooted maps, $k$-rooted 2-connected
maps and $k$-rooted 3-connected maps counted according to the number
of edges of their quotient map are respectively $\beta$-rational, $\eta$-rational, and
$\gamma$-rational. The explicit expressions are given in Appendix~\ref{sec:enum_edges}.

For $k\geq 2$, the series of $k$-rooted maps, $k$-rooted 2-connected
maps and $k$-rooted 3-connected maps counted according to the number
of vertices and the number of faces of their quotient map are respectively
$(\beta_1,\beta_2)$-rational, $(\eta_1,\eta_2)$-rational and
$(\gamma_1,\gamma_2)$-rational. The explicit expressions are given in Appendix~\ref{sec:enum_vertices}.

In particular, all these series are algebraic.
\end{theorem}

Using the algebraicity of the series of $k$-rooted maps,
methods of computer algebra can be used to quickly extract
their initial coefficients, see~\cite{Sa92}. Using Burnside's formula~(\ref{eq:burnsi}) and its
version in two variables if counting is done according to the
number of vertices and faces,
the enumeration of unrooted maps can be performed very efficiently: using
Maple, several hundreds of initial coefficients are easily computed. 
As in~\cite{Wa03}, the complexity models used here is that of 
\emph{arithmetic operation}, where an arithmetic operation is either 
the addition of two large integers both of size $\mathcal{O}(N)$ bits, 
or it is the multiplication or division of a large integer of size 
$\mathcal{O}(N)$ bits with a ``small'' integer of size $\mathcal{O}(\log(N))$ 
bits. 

\begin{theorem}
\label{theo:second}
For the enumeration of unrooted 2-connected and unrooted 3-connected maps
with respect to the number of edges, a table of the $N$ first
coefficients can be computed in $\mathcal{O}(N\log (N))$ operations.

For the enumeration of unrooted 2-connected and unrooted 3-connected maps
with respect to the number of vertices and faces, a table of the first 
coefficients with indices $(i,j)$ verifying $i+j\leq N$ can be computed 
in $\mathcal{O}(N^2)$ operations.
%The arithmetical operations involved in the calculations are, as in~\cite{Wa03}, the multiplication of a `''large'' integer with $\mathcal{O}(N)$ digits and of a ``small'' integer with $\mathcal{O}(\log (N))$ digits.
\end{theorem}
 For unrooted 2-connected maps, the same complexity results 
were obtained by Liskovets and Walsh~\cite{LW83}, with the difference that they give explicit 
formulas for the coefficients whereas we give explicit formulas for the 
generating functions. The improvement obtained by our method is for the 
family of 3-connected maps, which is
interesting as these objects correspond to oriented polyhedra. Walsh 
found no closed formula for the coefficients counting these maps, and 
proceeded with a costly procedure of iterative extraction of the 
coefficients from their cycle index sum. With our method, we still find 
explicit (algebraic) expressions for all generating functions of $k$-rooted 
3-connected maps. Our
complexity, in $\mathcal{O}(N\log (N))$ for one parameter and
$\mathcal{O}(N^2)$ for two parameters, improves significantly on the
complexity obtained by Walsh~\cite{Wa03}, of
$\mathcal{O}(N^3)$ for one parameter and
$\mathcal{O}(N^5)$ for two parameters.

\section{Definitions}
\label{sec:enum}
\subsection{Maps}
A map is a proper embedding of a connected graph (with possibly loops
and multiple edges) on a closed oriented surface, where \emph{proper}
means that edges are smooth arcs that do not cross. All the maps that are 
considered
in this article are on the sphere. For enumeration, maps
are considered up to orientation-preserving homeomorphisms of the topological 
sphere. Equivalently, two maps are identified if it is possible to obtain the 
second one from the first one by performing a continuous deformation of the 
sphere.
A map is said \emph{2-connected} (or non-separable) if it has no bi-partition of its edges intersecting at a single vertex. Equivalently, the two maps with one edge (i.e. the loop-map and the link-map) are 2-connected; and a map with at least two edges is 2-connected iff it has no loops and at
least two of its vertices have to be removed to disconnect the map. A map is 
said \emph{3-connected} if it has no loops nor
multiple edges and at
least three of its vertices have to be removed to disconnect the map.
A map is \emph{rooted} by marking and orienting one of its edges. This
operation suffices to eliminate all non trivial homeomorphisms of the
map. Hence, counting rooted maps is easier than counting maps because the
root can be used to start a recursive decomposition.
For $k\geq 2$, a \emph{$k$-rooted map} is a map with $k$ \emph{undistinguishable}
roots. This means that the $k$ objects obtained by marking differently
(say, in blue) one of the $k$ roots are equal. Rooted maps endowed with an
automorphism of order $k\geq 2$ are in bijection with $k$-rooted maps
(see~\cite{Li81} for more details). As $k$-rooted maps are easier to
handle for our purpose, we will manipulate them rather than rooted
maps with an automorphism of order $k$. 

\subsection{Quadrangulations}
A \emph{quadrangulation} is a map whose faces have degree 4.
A quadrangulation is said \emph{simple} if it has no multiple edge.
A quadrangulation is said \emph{irreducible} if each 4-cycle of edges
of the quadrangulation is the contour of one of its faces.
For each quadrangulation, its vertices can be colored in black and
white so that each edge connects a black and a white vertex. Such a
bicoloration is unique up to the choice of the colors. A
quadrangulation endowed with such a bicoloration is said
\emph{bicolored}.
A \emph{bi-rooted} quadrangulation is a rooted quadrangulation having
a secondary root. This secondary root is differently marked (say in blue), 
and is authorized to be equal to the primary root.  Taking the bicoloration 
into
account, a bi-rooted quadrangulation $Q$ is said
\emph{bicolor-consistent} if the origins of the two roots have the
same color when $Q$ is bicolored.

\subsection{Structure of $k$-rooted maps and method of quotient maps}
\label{section:krooted}

It was observed by Liskovets~\cite{Li81} that a $k$-rooted map can be
realized as an 
embedding on the geometrical sphere such that the embedding is invariant by a 
rotation of
angle $2\pi/k$ of the sphere \footnote{This point of view is not topologically
relevant but it helps to have a geometrical intuition and it gives an handy 
way to
define the quotient of a $k$-rooted map.}.  In addition, the two
points of intersection of the sphere with the rotation-axis are either a 
vertex or the centre
of a face if $k>2$, and can also be the middle of an edge if $k=2$. These
points are called the \emph{poles} of the $k$-rooted map. The
\emph{type} of a $k$-rooted map is the type of its two poles. For
example, if the two poles are a vertex and a face, then the $k$-rooted
map is said to have type face-vertex.

If we cut the sphere of the symmetrical embedding along two meridians forming a
dihedral angle of $2\pi /k$, we can extract a sector of the
map bordered by these two meridians. By
pasting together the two meridians, the sector becomes a map on the
sphere. The
symmetry of order $k$ of the initial geometrical embedding ensures
that this map is independant of the choice of the two meridians. We
call it the \emph{quotient-map} of the $k$-rooted map. Observe
that the quotient map has one root and two marked cells (the poles of
the $k$-rooted map). The method of quotient maps developed by
Liskovets consists in counting $k$-rooted maps of a family by studying
the structure of their quotient map. In the case of unconstrained
maps, it works well as quotient maps are essentially rooted maps
with two marked cells.

\subsection{Burnside's formula adapted to unrooted maps}
Consider a family of maps on the sphere (for example the family of
2-connected maps).
Let $c_n$, $c_n'$ and $c_n^{(k)}$ denote respectively the number of
unrooted, rooted and $k$-rooted ($k\geq 2$) maps of the family with $n$ edges.
Let $c_{ij}$, $c_{ij}'$ and $c_{ij}^{(k)}$ denote respectively the number of
unrooted, rooted and $k$-rooted ($k\geq 2$) maps of the family with $i+1$ vertices
and $j+1$ faces. Burnside's (orbit counting) formula was adapted by
Liskovets~\cite{Li81} to give the two following enumerative formulas for unrooted
maps, where $\phi(.)$ is Euler totient function.

\begin{equation}
\label{equa:burns}
2nc_n=c_n'+\sum_{k=2}^n \phi(k)c_n^{(k)}\ \ \ \ \ \ \ \ \ \ \ \ \  2(i+j)c_{ij}=c_{ij}'+\sum_{k=2}^{i+j}
\phi(k)c_{ij}^{(k)}
\end{equation}

As a consequence, counting unrooted maps in one parameter
(resp. two parameters) comes down to counting rooted maps
(already done for 2-connected and 3-connected maps, see~\cite{MS68}) and $k$-rooted
maps of the family with one parameter (resp. two parameters). 

\subsection{Bijection between maps and quadrangulations}
\label{sec:bijkrooted}
A classical result in map theory is a bijection between maps and
bicolored quadrangulations, that we shall refer to as the \emph{angular bijection}. We just detail its properties here. The angular bijection
is a bijection between maps with $n$ edges (resp. with $i$ vertices and $j$ faces) and
bicolored quadrangulations with $n$ faces (resp. with $i$ black and
$j$ white vertices). Indeed, by this bijection, vertices, faces and edges of
a map correspond respectively to black vertices, white vertices and
faces of the bicolored quadrangulation. 

In addition, under the angular bijection, rooted
maps are in bijection with rooted quadrangulations and $k$-rooted
maps are in bijection with so called 
\emph{$k$-rooted bicolored quadrangulations}, which
are defined as $k$-rooted quadrangulations such that the origins
of the $k$ roots have the same color when the quadrangulation is
bicolored. We will only deal with such $k$-rooted quadrangulations and
will shortly call them $k$-rooted quadrangulations. Observe that the
type of a $k$-rooted map and the type of its associated $k$-rooted
quadrangulation are linked by the above mentioned correspondence: for
example 2-rooted maps with type edge-face are in bijection with
2-rooted quadrangulations with type face-white vertex. We have 
seen in Section~\ref{section:krooted} that the two poles of a
$k$-rooted map are a face or a vertex if $k>2$ and can also be an edge
if $k=2$. Hence, a $k$-rooted quadrangulation can only have type 
vertex-vertex if
$k>2$, and can also have type face-face and type face-vertex if $k=2$.

Moreover, the angular bijection has the nice property that 2-connected
maps are in bijection with bicolored simple quadrangulations
and 3-connected maps are in bijection with bicolored irreducible 
quadrangulations. As a consequence, counting $k$-rooted 2-connected maps
according to the number of edges (resp. according to the numbers of vertices 
and faces) comes down to
counting $k$-rooted simple quadrangulations according to the number of
faces (resp. according to the numbers of black vertices and white vertices). 
The situation is the
same for 3-connected maps, but  with irreducible quadrangulations instead of
simple quadrangulations, see Figure~\ref{figure:diagschema}.

\subsection{Notations and conventions for the generating functions}
\label{sec:conv}
We will use the letters $F$, $G$ and $H$ to denote respectively
generating functions of $k$-rooted, $k$-rooted simple and $k$-rooted
irreducible quadrangulations. We will use the
subscripts $f$, $v$, $b$, and $w$ to denote respectively a pole which is a 
face, a vertex, a
black vertex and a white vertex. The subscripts $b$ and
$w$ are only
used for generating functions with two variables, where we have
to take the bicoloration into account. Moreover we will use the
exponent $(k)$ to denote a $k$-rooted quadrangulation. For example, 
$G_{vv}^{(k)}(y)$
is the series counting $k$-rooted simple
quadrangulations of type vertex-vertex by the
number of faces in their quotient map, and $H_{bw}^{(k)}(\zb,\zw)$ is
the series counting $k$-rooted
irreducible quadrangulations, whose poles are a black and a white vertex, 
by the number of black and white vertices in their quotient
map (and without counting the two axial vertices). 

More precisely, we adopt the following conventions for the generating
functions. These conventions will always be adopted for the series
counting rooted or $k$-rooted quadrangulations. The reader is advised to 
examine
them carefully only for the sake of verifying the later obtained
equations (Section~\ref{sec:count2connected} and 
Section~\ref{sec:count3connected}). For the series in one
variable, rooted quadrangulations will be counted according to their
number of faces, and $k$-rooted quadrangulations will be counted
according to the number of faces in their quotient map without counting the 
axial
faces. For example a $2$-rooted quadrangulation of type face-face with
$2n+2$ faces will
be counted in the coefficient of $x^n$ in the series $F_{f\!f}(x)$. For the
series in two variables, rooted quadrangulations will be counted
according to their number of black vertices minus 1 and to their
number of white vertices minus 1; $k$-rooted quadrangulations of type 
vertex-vertex
will be counted according to the number of black and white vertices in their 
quotient
map and without counting the two axial vertices; $k$-rooted
quadrangulations of type face-vertex will be counted according to the
number of  black and white vertices in their quotient map and without
counting the vertices incident to the axial face; $k$-rooted
quadrangulations of type face-face will be counted according to the
number of black vertices and white vertices of their quotient map and
without counting the vertices of one axial face. For example, a
$2$-rooted irreducible quadrangulation of type ``face-face'' with $2i+2$ black
vertices and $2j+2$ white vertices will be counted in the coefficient
of $\yb ^i\yw ^j$ in
the series $H_{f\!f}(\yb,\yw)$; a
$2$-rooted simple quadrangulation of type ``black vertex-face'' with $2i+1$ 
black
vertices and $2j+2$ white vertices will be counted in the coefficient
of $\yb ^i\yw ^j$ in
the series $G_{bf}(\yb,\yw)$; a $k$-rooted irreducible quadrangulation of type
``black vertex-white vertex'' with $ki+1$ black vertices and $kj+1$ white
vertices will be counted in the coefficient of $\zb ^i \zw^j$ in the
series $H_{bw}^{(k)}(\zb,\zw)$. 

Our conventions and the property that a quadrangulation with $n$ faces has $n+2$ vertices yield the following pleasant property: for each family of $k$-rooted maps, its generating function $f(w)$ in one variable and its generating function $f(\tb,\tw)$ in two variables are related by 
$$
f(w)=f(\tb,\tw).
$$ 
For example, $F_{f\!f}(x)=F_{f\!f}(x,x)$ and $H_{vv}(z)=H_{bb}(z,z)+H_{bw}(z,z)+H_{ww}(z,z)$.
 
%Moreover, given a series $f(\tb,\tw)$, we write $^{t}\!f$ for the series
%$^tf(\tb,\tw)=f(\tw,\tb)$.

\subsection{Algebraic structure of unconstrained $k$-rooted maps}
\label{sec:algunconstr}
\begin{lemma}
\label{lemma:betarational}
For $k\geq 2$, the generating functions of $k$-rooted quadrangulations in one
(resp. two) variable are $\beta$-rational
(resp. $(\beta_1,\beta_2)$-rational).
\end{lemma}
\begin{proof}
We take here only the example of $k$-rooted quadrangulations of type vertex-vertex (the treatment is similar but a little more involved when there is an axial face). From the method of quotient-map of Liskovets, the quotient of a
$k$-rooted quadrangulation of type vertex-vertex with $kn$ faces is a rooted
quadrangulation with $n$ faces and two marked vertices. A rooted 
quadrangulation with $n$ 
faces has $n+2$ vertices according to Euler relation, so that there are 
$(n+2)(n+1)/2$ 
possible choices of two poles. Hence, writing $F_n$ for the number of rooted 
quadrangulations with $n$ faces and $F_{vv,n}^{(k)}$ for the number of 
$k$-rooted 
quadrangulations of type vertex-vertex with $kn$ faces, we have 
$F_{vv,n}^{(k)}=\frac{1}{2}(n+2)(n+1)F_n$. Observe that it implies that
the series counting $k$-rooted quadrangulations of type vertex-vertex
by the number of faces of their quotient map does not depend on $k$. More 
precisely, 
we have $F_{vv}^{(k)}(x)=\frac{1}{2}x^2d^2F/dx^2+2xdF/dx+F$. The other series 
counting 
$k$-rooted quadrangulations also involve the first and second derivatives 
(or partial
derivatives for two variables) of the series $F$
counting rooted quadrangulations. This series is well-known to be
$\beta$-rational in one variable~\cite{BT} and $(\beta_1,\beta_2)$-rational 
in
two variables~\cite{Sc97}: $F(x)=\frac{\beta(2-9\beta)}{(1-3\beta)^2}$ and 
$F(\xb,\xw)=\frac{\beta_1+\beta_2-5\beta_1\beta_2+2\beta_1^2+2\beta_2^2}{(1-\beta_1-2\beta_2)(1-\beta_2-2\beta_1)}$. In addition, the property of being 
$\beta$-rational
(resp. $(\beta_1,\beta_2)$-rational) is stable
under derivation. Indeed, $dF/dx=(dF/d\beta)/(dx/d\beta)$ is the
quotient of two $\beta$-rational expressions, and we can proceed
similarly with two variables. The result follows.
\end{proof} 

\section{Tree-decompositions}
\label{section:treedecomp}
\subsection{Tree-decomposition by multiple edges}
\label{sec:treedecomp2}
\begin{figure}
\begin{center}
\includegraphics[width=13cm]{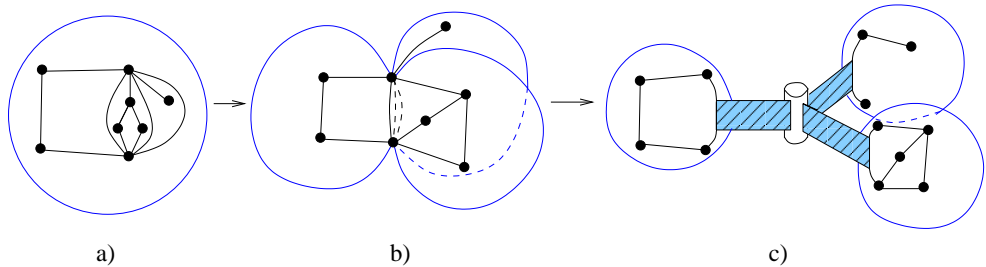}
\end{center}
\caption{Tree-decomposition by multiple edges of a
quadrangulation.}
\label{figure:decquad}
\end{figure}
We explain here how to transform an unrooted quadrangulation $Q$ (that
may have multiple edges) into a tree with two kinds of nodes: nodes
representing multiple edges and nodes representing simple
quadrangulations.

One way to see this decomposition is as follows. Take a multiple edge of
$Q$ of multiplicity $d$. Cut the sphere along each of the $d$ edges forming 
the multiple edge. In this way we obtain $d$ sectors, each
sector being delimited by two consecutive edges of the multiple
edge. Now, for each sector, identify the two meridians corresponding to the two
edges delimiting the sector by pasting them together. Thus we make out
of each sector a map on the sphere, and we can link these $d$ maps, at
their edge corresponding to the initial multiple edge, around a new
node: this will be the node of the tree corresponding to the multiple
edge. Now we can carry on recursively the tree-decomposition for each
of the $d$ maps, until all multiple edges have been split into nodes of
the tree.

Another way to see this decomposition is to imagine that we ``blow''
equally, from the interior of the sphere, each of the $d$ sectors
delimited by the multiple edge. We obtain thus $d$ components drawn
each on a sphere, where the $d$ spheres are connected (glued) at the
multiple edge, see Figure~\ref{figure:decquad}b. We can then represent
this multiple edge as a rigid link (see Figure~\ref{figure:decquad}c) around 
which the $d$ components are linked via
their unique edge belonging to the multiple edge. We can then
carry on the decomposition for each of the $d$ components.

\subsection{Tree-decomposition by separating 4-cycles}
\label{section:treedecompsecond}
The second tree-decomposition which we use consists in transforming a simple 
quadrangulation with at least 3
faces into a tree with two kinds of nodes: so-called \emph{axis-nodes} and
nodes corresponding to irreducible quadrangulations. The description
of this tree-decomposition can also be found in~\cite{Ku}. We describe
first in a recursive way the tree-decomposition for rooted objects, 
giving rise to a (rooted) decomposition-tree. Then, similarly as for 
the first tree-decomposition, we give a topological argument ensuring 
that the tree-decomposition can be equivalently performed on unrooted objects.

\begin{figure}
\begin{center}
\includegraphics[width=10cm]{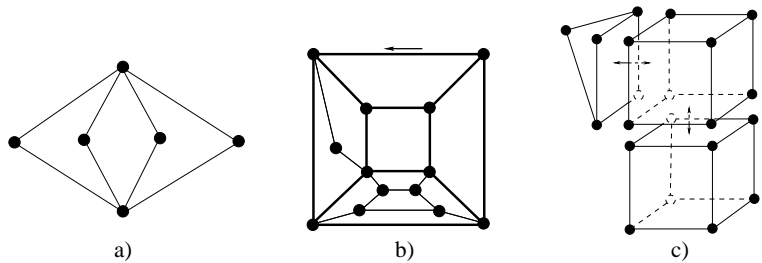}
\end{center}
\caption{An axis-map with 4 faces (a). The tree-decomposition of a 
quadrangulation by separating 4-cycles,
performed with a root (b) or without a root (c).}
\label{figure:decomp4cycle}
\end{figure}

For $k\geq 3$, we define the \emph{axis-map} with $k$ faces as the simple
quadrangulation consisting of two vertices
linked by $k$ parallel chains of 2 edges, each pair of two
consecutive paths forming one of the $k$ faces of the axis-map, see 
Figure~\ref{figure:decomp4cycle}a. The two vertices linked by the $k$ chains are called the \emph{extremal vertices} of the axis-map.

Now we state the following lemma of decomposition of a rooted simple
quadrangulation $Q$ with at least 3 faces:

\begin{lemma}
\label{lemma:decompaxis}
There exists a unique rooted quadrangulation $Q_0$, with maximal
possible number $k+1$ of faces such that:
\begin{itemize}
\item
$Q_0$ is an axis-map or an irreducible quadrangulation.
\item
There are $k$  rooted simple quadrangulations $Q_1,\ldots , Q_k$ with
at least 2 faces (including the rooted one) such that $Q$ can be seen
as the quadrangulation $Q_0$ where each of the $k$ non root faces
$f_i$ of $Q_0$ is substituted in a canonical way by one of the $Q_i$, 
$1\leq i \leq k$, the
contour of $f_i$ being replaced by the contour of the root face of
$Q_i$.
\end{itemize}
\end{lemma}
\begin{proof}
If there exists an internal chain of length 2 between two opposite
vertices of the outer face of $Q$, take the sequence of all chains of length
2 (including the two outer ones) between these two vertices. Forgetting
all other edges, we get an axis-map. Hence $Q$ can be seen as this axis-map 
where each non root face is
substituted by a quadrangulation.

Otherwise, define a proper 4-cycle of $Q$ as a 4-cycle
different from the contour of the root face of $Q$. Here we have to
see $Q$ as drawn in the plane with its root face as infinite face, so
that we can distinguish interior and exterior. A proper 4-cycle is
said maximal if it is not strictly included in the interior of any
other proper 4-cycle. As $Q$ has no path of length 2 connecting two opposite 
vertices of its outer face, it can easily be shown (see~\cite{MS68}) that the
interiors of maximal proper 4-cycles partition the interior of the outer face 
of
$Q$. Let $Q_0$ be the rooted quadrangulation obtained from $Q$ by keeping the
contour of the root face and of the maximal proper 4-cycles of
$Q$. The quadrangulation $Q_0$ is irreducible by maximallity
of the 4-cycles of which we have kept the contour. Hence we are in the
case where $Q$ can be seen as a rooted irreducible quadrangulation
where each inner face is substituted by a rooted quadrangulation.
\end{proof}

The first (resp. second) case of Lemma~\ref{lemma:decompaxis} correspond to 
the case
where the root node of the decomposition-tree is an axis-node
(resp. an irreducible quadrangulation). For example,
on Figure~\ref{figure:decomp4cycle}b, the rooted quadrangulation can be seen 
as a (rooted) cube
where two faces are substituted by another cube and by an axis-map with 3 faces. 

\textbf{Remark:} We make the following distinction when the root node of the 
decomposition-tree is an
axis-node: if the parallels chains of length 2 are incident to the
origin of the root, the root node of the tree is called a \emph{vertical} 
axis-node,
otherwise it is called an \emph{horizontal} axis-node.

Now we can carry on the tree-decomposition for each rooted
quadrangulation $Q_i$ with $1\leq i\leq k$. Thus, we get finally a (rooted) 
decomposition-tree
with two types of nodes: axis nodes and nodes that are irreducible
quadrangulations. Observe that, if $Q_0$ and the root node of one of
the $Q_i$ are simultaneously axis-nodes, then they are stretched in
perpendicular directions by maximallity of the number of faces of
$Q_0$.

The preceding decomposition on rooted objects is such that,
as in Section~\ref{sec:treedecomp2}, we can ``blow'' from the interior of the 
sphere to
``sculpt'' the quadrangulation $Q$ into a tree with nodes that are
irreducible quadrangulations and nodes that are axis-nodes, these nodes being
connected (glued) at so-called \emph{interconnection-faces}, see
Figure~\ref{figure:decomp4cycle}c. Hence we can say that an unrooted simple 
quadrangulation
``is'' its tree-decomposition after a judicious deformation of the
sphere. Thus we see that the topological shape  of the decomposition-tree in 
the
space does not depend on the face of the quadrangulation where we
choose to place the root to start the tree-decomposition. Hence an unrooted 
simple quadrangulation gives rise to an unrooted decomposition-tree.

\subsection{Centre of a tree}
The \emph{centre} of a tree $T$ is defined in the following recursive
way. If $T$ is reduced to an edge or a node, then the centre
of $T$ is this edge (resp. this node). Otherwise, remove all leaves of $T$ to
obtain a (shrinked) tree $\widetilde{T}$. Then the centre of $T$ is
defined to be the centre of $\widetilde{T}$. 

The important point is that the definition does not require that $T$ is
rooted. Hence the centre is invariant under any symmetry of $T$.

\section{Using the tree-decomposition by multiple edges to enumerate
unrooted 2-connected maps}
\label{sec:count2connected}
\subsection{Repercussion of the symmetry on the decomposition-tree}
\label{sec:introdecomp}

\begin{figure}
  \begin{center}
    \includegraphics[width=13cm]{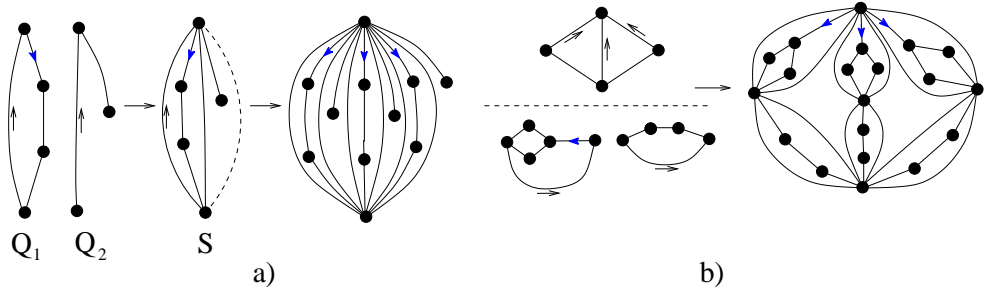}
  \end{center}
  \caption{Construction of a $k$-rooted quadrangulation of type $a$
(Figure a), and of type $b$ (Figure b).}
\label{figure:qqmult}
\end{figure}

As we have seen in Section~\ref{sec:treedecomp2}, the decomposition-tree of a
quadrangulation $Q$ is obtained by deforming the sphere in such a way that 
multiple edges can be split into link-nodes. This gives rise to a 
decomposition-tree ``living'' in the 3D-space. In addition, if
$Q$ is $k$-rooted ($k\geq 2$), then its decomposition-tree is invariant under 
the rotation-symmetry of order $k$ induced by its $k$-root.  

\begin{proposition}
\label{prop:centre}
For $k\geq 2$, the centre of the decomposition-tree of a $k$-rooted 
quadrangulation is a node (not an edge) of the tree, and it is the unique 
node of the tree fixed by the symmetry induced by the $k$-root. This node 
is called the \emph{core-node} of the decomposition-tree.
\end{proposition}
\begin{proof}
Assume that several nodes of the decomposition-tree are fixed by the 
symmetry. Then the axis of the rotation-symmetry has to pass by all these 
nodes. Hence, these nodes form a chain $x_0,\ldots,x_k$, where $x_i$ is 
connected to $x_{i+1}$ at a common edge $e_i$ of $Q$. As a consequence, 
the symmetry has to be of order 2 and to turn over such an edge $e_i$. This is 
impossible, because what we shortly call $k$-rooted quadrangulations are 
indeed $k$-rooted quadrangulations such that the origins of the $k$ roots 
have the same color when $Q$ is bicolored. In particular, the symmetry of 
the $k$-root can not exchange the two extremities of an edge of $Q$. Hence, 
at most one node of the tree can be fixed by the symmetry. Now assume that 
the centre of the decomposition-tree  is an edge $E$. This edge of the tree 
connects a node of type multiple-edge to a node that is a 
simple-quadrangulation $Q_s$ at an edge $e$ of $Q_s$. In addition, as 
noted in the definition of the centre of the tree, $E$ has to be fixed by 
the symmetry induced by the $k$-root. There are two ways the symmetry can 
fix $E$: either it exchanges its two extrimities, which is impossible as 
they are nodes of different types; or it fixes its two extremities, which 
is impossible because it would imply the presence of more than one node of 
the tree fixed by the symmetry.
\end{proof}
There are two possibilities for the core-node: either it is a node of type 
multiple edge -we say that $Q$ has type $a$- or it is a node of type simple 
quadrangulation -we say that $Q$ has
type $b$-.

\subsection{Construction of a $k$-rooted quadrangulation of type $a$.}
\label{sec:quad_a}
First we need to define a \emph{simply rooted} quadrangulation as a
quadrangulation whose root edge does not belong to a multiple edge.
% A bi-rooted
% quadrangulation is said \emph{bicolor-consistent} if the origins of
% its two roots have the same color when the quadrangulation is
% bicolored.

Now we explain how to construct a $k$-rooted quadrangulation whose
centre of the decomposition-tree is a multiple edge with multiplicity
$k\cdot d$ ($d \geq 1$), see Figure~\ref{figure:qqmult}a. Take a bi-rooted, 
simply-rooted (i.e. whose primary root is a simple edge)
quadrangulation $Q_1$. Cut it along its primary root-edge, thus
transforming $Q_1$ into a sector with two bording meridians. For
convenience, we consider the root as an arrow placed slightly on the
right of the root edge, so that cutting along the root edge does not
``split'' the root.  Among
these two meridians, we call root-meridian the one having the root
slightly on its right after the cutting. 

Now take $d-1$ simply rooted
quadrangulations $Q_2,\ldots,Q_d$ and perform the same cutting
operation on them as on $Q_1$. Then paste the root meridian of $Q_2$ with the
non-root meridian of $Q_1$, the pasting operation being such that the
orientations of the roots of the two sectors coincide. Then,
iteratively for each $2\leq i\leq d$, paste the root meridian of $Q_i$ with
the non-root meridian of $Q_{i-1}$, so as to obtain a big sector $S$
whose root meridian is the root meridian of $Q_1$. Now make $k$ copies
$S_1,\ldots ,S_k$ of $S$ and, for each $1\leq i\leq k$, paste the root
meridian of $S_i$ with the non-root meridian of $S_{i-1}$. In this way
we obtain finally a quadrangulation (on the sphere) consisting of $k$ identical
sectors, each carrying a blue root (the secondary root of $Q_1$). By
erasing the mark of the primary root of $Q_1$ and of the roots of
$Q_2\ldots Q_d$ in each sector, we obtain a $k$-rooted quadrangulation
of type $a$, see Figure~\ref{figure:qqmult}a. 
%In addition, it is easy to see that this construction is
%inversible, so that  each $k$-rooted quadrangulation of type $a$ is
%obtained in a unique way by such a construction.
Observe that each $k$-rooted quadrangulation of type $a$ is obtained
exactly twice by this construction. Indeed, the inverse operation
consists in choosing an extremity $v$ (two possibilities) of the
central multiple edge  and then orienting all edges of the multiple
edge toward $v$. 

Observe also that, if $Q_1$ is taken to be bi-rooted bicolor-consistent, 
then the 2-fold ambiguity of the construction
disappears and it becomes a bijection.

We write $f(x)$ for the series counting simply rooted quadrangulations by
their number of faces. Observe that the generating function $F(x)$ of rooted 
quadrangulations and $f(x)$ are related by $F(x)=f(x)/(1-f(x))$, following 
from the fact that a rooted quadrangulation whose root edge has mutliplicity 
$d$ can be decomposed in $d$ simply rooted quadrangulations. The construction 
given above allows us to express the generating function $F_{vv}^{(k)a}(x)$ of
$k$-rooted quadrangulations of type $a$ in terms of $f(x)$:

$$
F_{vv}^{(k)a}(x)=\frac{1}{2}(4xf'(x)).\frac{1}{1-f(x)},
$$ 
where we use the subscript $vv$ for $F_{vv}^{(k)a}(x)$ because all $k$-rooted 
quadrangulations of type $a$ clearly have type vertex-vertex.

\subsection{Construction of a $k$-rooted quadrangulation of type $b$.}
\label{sec:quad_b}
In this section, we give a construction of a $k$-rooted quadrangulation whose core-node is a simple quadrangulation as a composed object, see Figure~\ref{figure:qqmult}b. Take a $k$-rooted 
simple quadrangulation
$Q_s$. For the $k$-orbite of root edges, either leave its $k$ edges
untouched (Case 1) or perform the following operation (Case 2): take a 
bi-rooted bicolor-consistent quadrangulation $\widetilde{Q}$. Then cut $Q_s$
along each of its $k$ root edges and cut $\widetilde{Q}$ along its
primary root edge, transforming $\widetilde{Q}$ into a sector $S$ bordered by 
two
meridians. Take $k$ copies of $S$ and for each (cutted)
root-edge $e$ of $Q_s$, place a copy of $S$ in the empty sector
of $Q_s$ 
leaved by the cutting of $e$. This substitution is done by pasting the two 
meridians of
$S$ with the two border-edges of $Q_s$ created by cutting
 $e$, and by making the orientation of $e$ and of the primary root edge of
$S$ coincide.

Proceed similarly for each $k$-orbite of non-root edges of $Q_s$, with
the only difference that the quadrangulation used for
the substitution is not bi-rooted but just rooted. Finally, keep only the $k$ 
marks of the roots of $Q_s$ if we are in Case 1 (i.e. no substitution at the 
root edges of $Q_s$), and keep only the 
marks of the secondary roots of the $k$ copies of $\widetilde{Q}$ if we are 
in Case 2. Thus, we
obtain a $k$-rooted quadrangulation $Q$ of type $b$. 
%It is easy to see
%that this construction is invertible, i.e. that we can guess the roots
%of $Q$ whose marks we have erased.

Observe that $k$-rooted quadrangulations of type $b$ obtained by this
construction always have the following property: their $k$ root edges are 
simple if their
incident face (the face on their right) belongs to the central simple
quadrangulation (because this case corresponds to Case 1 where there is no 
substitution at the root edges of $Q_s$). The missing $k$-rooted 
quadrangulations of type $b$, i.e. those whose root edges are not simple and 
are incident to the central simple quadrangulation,
are obtained by the same construction, with the difference that we
always cut the $k$ root edges of $Q_s$. Then the other difference is that the 
first substituted
quadrangulation $\widetilde{Q}$ is not bi-rooted but just rooted. At
the end of this construction, we only keep the mark of the $k$ roots
of $Q_s$.

Similarly as in Section~\ref{sec:quad_a}, these two complementary
constructions allow us to obtain all
$k$-rooted quadrangulations of type $b$ in a bijective way. We introduce the family 
$\mathcal{E}$ as the union of the set of bi-rooted bicolor consistent
quadrangulations, of the set of rooted
quadrangulations, and of the
link-map (the map with one edge). The construction can be summarized by saying
that the root edges of $Q_s$ are substituted by $k$ copies of an
object of $\mathcal{E}$. We write $E(x)$ for the series counting
objects of $\mathcal{E}$ by
their number of faces, so that $E(x)=2xF'(x)+F(x)+1$.

We write $G_{vv,n}^{(k)}$ for the number of
 $k$-rooted simple quadrangulations with $n$ $k$-orbites of
faces. A quadrangulation has always doubly more edges
than faces, so that an object counted by $G_{vv,n}^{(k)}$ has $2n$
orbites of $k$ edges.    

The construction by substitution gives rise to the following equation for the 
case
where $Q_s$ has type vertex-vertex:
$$F_{vv}^{(k)b}(x):=\sum_n F_{vv,n}^{(k)b}x^n=\sum_n G_{vv,n}^{(k)}E(x)\left( 1+F(x)
\right) ^{2n-1}=\frac{E(x)}{1+F(x)}G_{vv}^{(k)}\left((1+F(x))^2\right)$$

Similarly, the two following expressions can be obtained,
corresponding respectively to the case where $Q_s$ has type
face-vertex and type  face-face:
  
\begin{minipage}{6cm}
$$E(x)G_{vf}^{(2)}\left((1+F(x))^2\right) ,$$
\end{minipage}  
\begin{minipage}{6cm}
$$E(x)(1+F(x))G_{f\!f}^{(2)}\left((1+F(x))^2\right) .$$
\end{minipage}  

\subsection{Obtaining the equations}
As $k$-rooted quadrangulations are partitioned into two sets whether their core-node is a multiple edge or a simple
quadrangulation,  we obtain the following equations by taking the sum of the 
series obtained in Section~\ref{sec:quad_a}
and Section~\ref{sec:quad_b}:

\begin{eqnarray}
\label{eq:2connected}
F_{vv}^{(k)}(x)&=&2\frac{xf'(x)}{1-f(x)}+\frac{E(x)}{1+F(x)}G_{vv}^{(k)}\left((1+F(x))^2\right)\label{eq:2connectedvv}\\
F_{vf}^{(2)}(x)&=&E(x)G_{vf}^{(2)}\left((1+F(x))^2\right)\label{eq:2connectedfv}\\
F_{f\!f}^{(2)}(x)&=&E(x)(1+F(x))G_{f\!f}^{(2)}\left((1+F(x))^2\right),\label{eq:2connected_ff}
\end{eqnarray}
where the only unknown series are $G_{vv}^{(k)}$, $G_{vf}^{(2)}$ and 
$G_{f\!f}^{(2)}$.

Similar equations can be easily  obtained in two variables by taking
the bicoloration  of vertices into account. We define
$df(\xb,\xw)=\frac{d}{dt}f(t\xb,t\xw)_{t=1}$ and write $E(\xb,\xw)$
for the series in two variables counting the family $\mathcal{E}$. It
is easy to establish, having conventions of Section~\ref{sec:conv} in mind and
using Euler relation (i.e. the sum of vertices and faces is the number
of edges +2), that 
$E(\xb,\xw)=2\frac{d}{dt}F(t\xb,t\xw)_{t=1}+F(\xb,\xw)+1$. Then 
Equation~(\ref{eq:2connectedvv}) becomes for example:
\begin{equation}
\left\{
\begin{array}{rcl}
F_{bw}^{(k)}(\xb,\xw)&=&2\frac{df}{1-f}+\frac{E}{1+F}G_{bw}^{(k)}\left(\xb(1+F)^2,\xw(1+F)^2\right)\label{eq:2connected2variables}\\
F_{bb}^{(k)}(\xb,\xw)&=&\frac{E}{1+F}G_{bb}^{(k)}\left(\xb(1+F)^2,\xw(1+F)^2\right)\\
F_{ww}^{(k)}(\xb,\xw)&=&\frac{E}{1+F} G_{ww}^{(k)}\left(\xb(1+F)^2,\xw(1+F)^2\right),
\end{array}
\right.
\end{equation}
where all series (including $f$ and $F$) have two variables, one for
the number of black vertices, the other one for the number of white
vertices.

As observed in Lemma~\ref{lemma:betarational}, the series $F_{vv}^{(k)}$ 
(and it is also the case for $F_{bw}^{(k)}$, $F_{bb}^{(k)}$ and 
$F_{ww}^{(k)}$)  does not
depend on $k$. Hence it follows from the form of Equation~(\ref{eq:2connectedvv}) and~(\ref{eq:2connected2variables}) that the series $G_{vv}^{(k)}$ 
(and also the series $G_{bw}^{(k)}$, $G_{bb}^{(k)}$ and $G_{ww}^{(k)}$) 
does not depend on $k$. Hence the
exponent $(k)$ can be ommited.

\begin{lemma}
\label{lemma:etarational}
For $k\geq 2$, the generating functions
of $k$-rooted simple quadrangulations in one variable (resp. two
variables) are $\eta$-rational (resp. $(\eta_1,\eta_2)$-rational).
\end{lemma}
\begin{proof}
From Lemma~\ref{lemma:betarational}, we know that $F_{vv}(x)$, $F_{vf}(x)$ and $F_{f\!f}(x)$ are $\beta$-rational, and so
are $x$ (because $x=\beta-3\beta^2$), $F(x)$ (as proved in~\cite{Tu63}), $f(x)$ (because $F=f/(1-f)$), 
and $E(x)$ (because $\beta$-rationality is stable under derivation). Hence it follows from
Equations~(\ref{eq:2connectedvv}),(~\ref{eq:2connectedfv}) 
and~(\ref{eq:2connected_ff}) that $G_{vv}\left(x(1+F)^2\right)$, 
$G_{vf}\left( x(1+F)^2\right)$ and $G_{f\!f}\left(
x(1+F)^2\right)$ are $\beta$-rational. Now we have to make the change
of variable $y=x(1+F)^2$. It was observed in~\cite{BT} that 
$\beta(x)=\eta(y)/(1+3\eta(y))$ when $y$ and $x$
are linked by the change of variable
$y=x(1+F)^2$. Hence, replacing $\beta(x)$ by
$\eta(y)/(1+3\eta(y))$ in the respective $\beta$-rational expression of 
$G_{vv}\left(
x(1+F)^2\right)$, $G_{vf}\left( x(1+F)^2\right)$ and $G_{\!f\!f}\left(
x(1+F)^2\right)$, we obtain $\eta$-rational expressions for
$G_{vv}(y)$, $G_{vf}(y)$ and $G_{f\!f}(y)$.

We can proceed similarly in two variables, using the fact that
$\beta_1(\xb,\xw)$ and $\beta_2(\xb,\xw)$ have a rational expression
in terms of $\eta_1(\yb,\yw)$ and $\eta_2(\yb,\yw)$ when $(\yb,\yw)$
and $(\xb,\xw)$ are linked by the change of variable $(\yb,\yw)=(\xb
(1+F)^2,\xw (1+F)^2)$.
\end{proof}

\begin{lemma}
\label{lemma:complexity2}
The $N$ initial coefficients counting unrooted 2-connected maps
according to their number of edges can be computed with
$\mathcal{O}(N\log (N))$ operations.

The table of initial coefficients with indices $(i,j)$ and $i+j\leq N$
counting unrooted 2-connected maps according to their number of vertices and
faces can be computed with $\mathcal{O}(N^2)$ operations.
\end{lemma}

\begin{proof}
First we use the following notation. For a series $f$ in one variable
(resp. two variables), we denote by $\mathcal{C}_N(f)$ the number of
operations necessary to compute its $N$ initial coefficients
(resp. its coefficients with indices $(i,j)$ and $i+j\leq N$).  
Writing $g_n$ (resp. $g_{ij}$) for the number of unrooted 2-connected maps
with $n$ edges (resp. $i+1$ vertices and $j+1$ faces),
Burnside's formula~(\ref{equa:burns}) can easily be transposed
in the following equations on series:

\begin{eqnarray*}
\sum_n 2ng_ny^n&=&G(y)+yG_{vf}(y^2)+y^2G_{f\!f}(y^2)+\sum_{k\geq
2}\phi(k)G_{vv}(y^k)\\
\sum_{i,j} 2(i+j)g_{ij}\yb ^i\yw ^j&=&G(\yb,\yw)+\yw G_{bf}(\yb ^2,\yw
^2)+\yb G_{wf}(\yb ^2,\yw ^2)+\yb\yw G_{f\!f}(\yb ^2,\yw ^2)\\
&&+\sum_{k\geq
2}\phi(k)\left( \frac{\yb}{\yw} G_{bb}(\yb ^k,\yw ^k)+G_{bw}(\yb ^k,\yw ^k)+\frac{\yw}{\yb}
G_{ww}(\yb ^k,\yw ^k)\right)
\end{eqnarray*}

According to~\cite{BT}, $G(y)$ is $\eta$-rational; and according to Lemma~\ref{lemma:etarational}, $G_{vf}(y)$, $G_{f\!f}(y)$ and
$G_{vv}(y)$ are $\eta$-rational. Hence these series are algebraic  because
they live in the algebraic extension of the algebraic series
$\eta(y)$. Hence, they are differentiably finite
(see~\cite{Sta}), i.e. solution of a linear differential equation with
polynomial coefficients. Taking coefficient $[y^n]$ in this
differential equation yields a linear recurrence with polynomial 
coefficients for the coefficients of these series.  As a
consequence, the $N$ initial coefficients of these series can be
computed with $\mathcal{O}(N)$ ``arithmetical'' operations. 
Hence, $\mathcal{C}_N\left( \sum 2ng_n
\right) =\mathcal{C}_{N}(G)+\mathcal{C}_{N/2}(G_{vf}+G_{f\!f})+\sum_{k=2}^N\mathcal{C}_{N/k}(G_{vv})=\mathcal{O}(N)+\mathcal{O}(N/2)+\sum_{k=2}^N\mathcal{O}(N/k)=\mathcal{O}(N\log(N))$.

Similarly, an algebraic series in two variables is also D-finite. 
Hence its coefficients verify two linear recurrences, one for each
index. As a consequence, if $f(\yb,\yw)$ is algebraic, then
$\mathcal{C}_N(f)=\mathcal{O}(N^2)$. As the series of rooted and $k$-rooted simple
quadrangulations in two variables are $(\eta_1,\eta_2)$-rational, they are
algebraic. Hence, 
%\begin{eqnarray*}
%\mathcal{C}_N\left( \sum_{i,j}2(i+j)c_{ij}\right) &=&\mathcal{C}_N(G)+\mathcal{C}_{N/2}(G_{f\!f}+G_{bf}+G_{wf})+\sum_{k=2}^N\mathcal{C}_{N/k}(G_{bb}+G_{bw}+G_{ww})\\
%&=&\mathcal{O}(N)+\mathcal{O}((N/2)^2)+\sum_{k=2}^N\mathcal{O}((N/k)^2)=\mathcal{O}(N^2)
%\end{eqnarray*}
$\mathcal{C}_N\left( \sum_{i,j}2(i+j)g_{ij}\right) =\mathcal{C}_N(G)+\mathcal{C}_{N/2}(G_{\!f\!f}+G_{bf}+G_{wf})+\sum_{k=2}^N\mathcal{C}_{N/k}(G_{bb}+G_{bw}+G_{ww})=\mathcal{O}(N)+\mathcal{O}((N/2)^2)+\sum_{k=2}^N\mathcal{O}((N/k)^2)=\mathcal{O}(N^2)$,
where we use the fact that $\sum_k 1/k^2$ converges. 
\end{proof}

\section{Using the tree-decomposition by separating 4-cycles to enumerate
unrooted 3-connected maps}
\label{sec:count3connected}

\subsection{Introduction}
In this part, we use the tree-decomposition by separating 4-cycles explained in
Section~\ref{section:treedecompsecond}. 
This tree-decomposition states that a simple quadrangulation can be 
seen as an arborescent structure whose nodes are either 
irreducible quadrangulations or so-called axis-nodes. From this
decomposition, we will obtain 
equations linking generating functions of $k$-rooted irreducible
quadrangulations and 
generating functions of $k$-rooted simple quadrangulations. As we
have already obtained expressions for the generating functions of $k$-rooted 
simple quadrangulations 
in Section~\ref{sec:count2connected}, 
 we will obtain from these equations the generating functions of
$k$-rooted 
irreducible quadrangulations, from which unrooted
3-connected maps can be enumerated
using the angular bijection and
Burnside's 
formula, see Figure~\ref{figure:diagschema}.

We treat first the case of $k$-rooted objects with $k\geq 3$. 
The case of 2-rooted objects is more difficult
 (for example a symmetry of order 2 of an axis-map can exchange its
extremal vertices), and will be thoroughly treated in Section~\ref{sec:2rootedirreducible}.

\subsection{The case of $k$-rooted irreducible quadrangulations with 
$k\geq 3$.}
\label{section:krootedsimple}
A first important remark is that all $k$-rooted quadrangulations have type
 vertex-vertex for $k>2$, as we have seen in Section~\ref{section:krooted}.
We introduce the families $\mathcal{W}$ of rooted simple
quadrangulations with at least two faces (this excludes the degenerated 
one-face quadrangulation) and the family $\mathcal{J}$
consisting of the objects of $\mathcal{W}$ whose root node of the
decomposition tree is not an horizontal axis-node. We write $W(y)$ and
$J(y)$ for the series counting these two families by their number of
faces. Observe that $W(y)=G(y)-2y$, where $G(y)$ is the series of rooted 
simple quadrangulations, and
$W(y)/y=\frac{J(y)/y}{1-J(y)/y}$. We define also the families
$\mathcal{W}'$ and $\mathcal{J}'$ of objects of $\mathcal{W}$ and
$\mathcal{J}$ having a secondary root incident to a face different
from the root face. The series counting objects of $\mathcal{W}'$ and
$\mathcal{J}'$ by their number of faces are respectively $4C(y)$ and
$4B(y)$ where $C(y)=y\frac{d}{dy}W(y)-W(y)$ and 
$B(y)=y\frac{d}{dy}J(y)-J(y)$. 

Let $Q$ be a $k$-rooted simple quadrangulation ($k\geq 3$) with at least 3
faces. The decomposition-tree of $Q$ is
invariant under the symmetry of order $k$ induced by the $k$-root of
$Q$. As $k>2$, the centre of the decmposition-tree is a node (not an edge) 
and is the unique node  
invariant by the symmetry. We call this node the \emph{core-node} 
(we will see later that the definition of the core-node requires more 
attention for 2-rooted objects).  Two cases can arise: either the
core-node is an axis-node -we say that $Q$ has type $\mathfrak{a}$- or it is an
irreducible quadrangulation -we say that $Q$ has type $\mathfrak{b}$-.

\subsubsection{Construction of $k$-rooted simple quadrangulations of type
$\mathfrak{a}$}
\label{sec:quadsimp_a}

%\begin{figure}
%  \begin{center}
%    \includegraphics[width=13cm]{Figures/simpleaxis.eps}
%  \end{center}
%  \caption{Construction of a $k$-rooted simple quadrangulation of type $a$
%(Figure a), and of type $b$ (Figure b).}
%\label{figure:simpleaxis}
%\end{figure}

Similarly as in Section~\ref{sec:quad_a}, we give a construction, in terms 
of a composed object, of a $k$-rooted simple
quadrangulation whose core-node is an axis-map
with $k\cdot d$ faces. Take a $k$-rooted axis-map
with $k\cdot d$ faces and whose roots point toward the same extremal vertex of the
axis-map, which we call the \emph{pointed extremal vertex}. Then take $k$ copies of
an object $Q_1$ of $\mathcal{J}'$ and substitute each root face of the
axis-map by one of these copies, making the primary root of the copies
of $Q_1$ be oriented toward the pointed extremal vertex of the axis-map. Proceed
similarly for each $k$-orbite of non-root faces of the axis-map, with
the only difference that the substituted objects are $k$ copies of an
object of $\mathcal{J}$ instead of $\mathcal{J}'$. Finally keep only
the marks of the secondary root of the $k$ copies of $Q_1$.

As in Section~\ref{sec:quad_a},  each $k$-rooted simple
quadrangulation of type $\mathfrak{a}$ is obtained exactly twice by this
construction. Hence, the series counting $k$-rooted simple
quadrangulations of type $\mathfrak{a}$ is:

$$
G_{vv}^{(k)\mathfrak{a}}(y)=2\frac{B(y)}{y}\frac{1}{1-J(y)/y}.
$$

\subsubsection{Construction of $k$-rooted simple quadrangulations of type
$\mathfrak{b}$}
\label{sec:quadsimpb}
 As precedently, we give a construction of $k$-rooted simple
quadrangulations of type $\mathfrak{b}$ as composed objects. Take a $k$-rooted
irreducible quadrangulation $Q_{irr}$. Take $k$ copies of an object
$Q_1$ of $\mathcal{W}'$ and substitute each root face of $Q_{irr}$ by
one of the copies of $Q_1$ in a ``canonical'' way, e.g. by superposing
the primary root edge of $Q_1$ with the root edge of the face where the
substitution takes place. Then proceed similarly for each $k$-orbite
of non-root faces of $Q_{irr}$, with the difference that the substituted 
objects are
$k$ copies of an object of $\mathcal{W}$ instead of
$\mathcal{W}'$. Finally keep only the marks of the secondary root of
the $k$ copies of $Q_1$. 

By this construction, all $k$-rooted simple quadrangulations of type
$\mathfrak{b}$ are
obtained exactly 4 times. Indeed, as a quadrangular face has 4 sides,
there are 4 possibilities to guess the primary root edge of the $k$
copies of $Q_1$. Hence the series counting $k$-rooted simple 
quadrangulations of type $\mathfrak{b}$ is given by the following 
expression, involving the series $H_{vv}^{(k)}(z)$ of $k$-rooted 
irreducible quadrangulations:

$$
G_{vv}^{(k)\mathfrak{b}}=\frac{1}{4}\sum_n H_{vv,n}^{(k)}\frac{4C(y)}{y}\left(\frac{W(y)}{y}\right)^{n-1}=\frac{C(y)}{W(y)}H_{vv}^{(k)}(W(y)/y)
$$

\subsubsection{Obtaining the equations}

The set of $k$-rooted simple quadrangulations is partitioned into two sets
whether the core-node is an axis-node or an
irreducible quadrangulation. Hence, summing the series obtained in
Section~\ref{sec:quadsimp_a} and Section~\ref{sec:quadsimpb}, we
obtain the following equation linking the series of $k$-rooted simple
quadrangulations with the series of $k$-rooted irreducible
quadrangulations, for $k>2$:

\begin{equation}
\label{eq:kirr}
G_{vv}^{(k)}(y)=2\frac{B(y)}{y}\frac{1}{1-J(y)/y}+\frac{C(y)}{W(y)}H_{vv}^{(k)}(W(y)/y)
\end{equation}

Similar equations can easily be obtained in two variables by taking the 
bicoloration of $Q$ into account. Writing $C(\yb,\yw)=\yb\frac{\partial
W}{\partial \yb}+\yw\frac{\partial
W}{\partial \yw}-W$ and $B(\yb,\yw)=\yb\frac{\partial
J}{\partial \yb}+\yw\frac{\partial
J}{\partial \yw}-J$ for the versions in two variables of $C(y)$ and
$B(y)$,
the version in two variables of Equation~(\ref{eq:kirr}) is

\begin{eqnarray}
\label{eq:kirr2}
G_{bb}^{(k)}(\yb,\yw)&=&\frac{B}{\yb}\frac{1}{1-J/\yb}+\frac{C}{W}H_{bb}^{(k)}(W/\yw,W/\yb)\\
G_{ww}^{(k)}(\yb,\yw)&=&\frac{B}{\yw}\frac{1}{1-J/\yw}+\frac{C}{W}H_{ww}^{(k)}(W/\yw,W/\yb)\\
G_{bw}^{(k)}(\yb,\yw)&=&\frac{C}{W}H_{bw}^{(k)}(W/\yw,W/\yb).
\end{eqnarray}
Observe that these equations are the same for all values 
 of $k$. As we have already seen that $G_{vv}^{(k)}(y)$ does not
depend on $k$, $H_{vv}^{(k)}(z)$ does also not depend on $k$, so that
we can denote this series by $H_{vv}^{\geq 3}$. We can make the
same remark for the series in two variables and adopt the same
notation for the exponent. 
\begin{lemma}
\label{lemma:gammarational}
For $k\geq 3$, the series of $k$-rooted irreducible quadrangulations in one 
variable (resp. two variables) is $\gamma$-rational 
(resp. $(\gamma_1,\gamma_2)$-rational) and does not depend on $k$. 
\end{lemma}
\begin{proof}
The proof is similar as the proof of Lemma~\ref{lemma:etarational}. In one
variable, we use the form
of Equation~(\ref{eq:kirr}) to see that $H_{vv}^{\geq 3}(W(y)/y)$ is
$\eta$-rational. Indeed all series appearing in Equation~(\ref{eq:kirr}) and
different from $H_{vv}^{\geq 3}(W(y)/y)$ have an explicit 
$\eta$-rational expression, as these series involve the series $W$ (which is
$\eta$-rational from~\cite{MS68}) or the series $J$ (equal to $W/(1+W/y)$) or
their derivatives, which are also $\eta$-rational because
$dW/dy=(dW/d\eta)/(dy/d\eta)$. Then we use the fact, shown in~\cite{MS68}, 
that $\eta(y)=\gamma(z)/(2\gamma(z)+1)$ when $z$ and $y$ are linked 
by the
change of variable $z=W(y)/y$. Substituting $\eta$ by
$\gamma/(2\gamma+1)$ in the
$\eta$-rational expression of $H_{vv}^{\geq 3}(W(y)/y)$, we obtain a
$\gamma$-rational expression for $H_{vv}^{\geq 3}(z)$. 

The proof for two variables is similar, using in particular the fact
that $\eta_1(\yb,\yw)$ and $\eta_2(\yb,\yw)$ have a rational
expression in terms of $\gamma_1(\zb,\zw)$ and $\gamma_2(\zb,\zw)$
when $(\zb,\zw)$ and $(\yb,\yw)$ are linked by the change of variable 
$(\zb,\zw)=(W/\yw ,W/\yb )$.
\end{proof}

%\begin{lemma}
%\label{lemma:complexity3}
%The $N$ initial coefficients counting unrooted 3-connected maps
%according to their number of edges can be computed with
%$\mathcal{O}(N\log (N))$ operations.

%The table of initial coefficients with indices $(i,j)$ and $i+j\leq N$
%counting unrooted 3-connected maps according to their number of vertices and
%faces can be computed with $\mathcal{O}(N^2)$ operations.
%\end{lemma}
%\begin{proof}
%Using the algebraicity of the generating function of $k$-rooted
%irreducible quadrangulations, we can perform the same treatment as in
%the proof of Lemma~\ref{lemma:complexity2}.
%\end{proof}

%Finally, Lemma~\ref{lemma:complexity2} and~\ref{lemma:complexity3}
%yield Theorem~\ref{theo:second}. Using Tutte's bijection between
%$k$-rooted objects (see also Figure~\ref{figure:diagschema}), Lemma~\ref{lemma:betarational},~\ref{lemma:etarational} and~\ref{lemma:gammarational} 
%yield Theorem~\ref{theo:first}.

\subsection{The case of 2-rooted irreducible quadrangulations}
\label{sec:2rootedirreducible}
\subsubsection{Introduction}
The case of 2-rooted objects requires a careful treatment. As we have seen in
Section~\ref{section:krooted}, a
2-rooted quadrangulations can have one or two axial faces. In addition, as 
opposed to the case $k>2$, other nodes than the centre of
the decomposition tree can be fixed by the symmetry of order 2 induced
by the 2-root. Another difficulty is the fact that a symmetry of order
2 of an axis-map can exchange its extremal vertices.

The object of this section is to show a result similar to
Lemma~\ref{lemma:gammarational} for 2-rooted objects, i.e. we would like to 
show first that the series of 2-rooted irreducible
quadrangulations in one variable (resp. two variables), composed with
$W(y)/y$ (resp. $\yW$) are $\eta$-rational
(resp. $(\eta_1,\eta_2)$-rational). From that, as in 
Lemma~\ref{lemma:gammarational}, it will follow that all series of 2-rooted 
irreducible quadrangulations in one variable (resp. two variables) are 
$\gamma$-rational (resp. $(\gamma_1,\gamma_2)$-rational). In order to prove 
this result, we have to find simple equations linking the series of 2-rooted simple quadrangulations and the series of
2-rooted irreducible quadrangulations. This requires to find a convenient 
partition of 2-rooted irreducible and 2-rooted simple
quadrangulations in several families. 
%Then we have to find equations
%involving these series and ensuring that, similarly as for $k$-rooted
%objects with $k\geq 3$, the series of 2-rooted irreducible
%quadrangulations in one variable (resp. two variables), composed with
%$W(y)/y$ (resp. $\yW$) are $\eta$-rational
%(resp. $(\eta_1,\eta_2)$-rational). 

\subsubsection{Introduction of families of 2-rooted simple and irreducible 
quadrangulations.}
We define here 5 families of 2-rooted simple quadrangulations and 5
families of 2-rooted irreducible quadrangulations that partition
respectively the set of
2-rooted simple and the set of 2-rooted irreducible quadrangulations.

We will obtain 5 equations linking the 5 generating functions of
irreducible quadrangulations, 
which are unknown, and the 5 generating functions of 2-rooted simple
quadrangulations, 
which are known. Moreover, we will see that this system of 5 equations
is upper triangular, 
so that it is easy to solve.

We define:
\begin{itemize}
\item
The families $\Gvv$ (resp. $\Qvv$) of 2-rooted simple
(resp. irreducible) 
quadrangulations of type \emph{vertex-vertex}. We write $\gvv (y)$
(resp. $\qvv (z)$) for their generating functions. If the
bicoloration is taken into account, $\Gvv$ (resp. $\Qvv$) is partitioned in 
three 
families
$\Gbb$, $\Gww$, $\Gbw$ (resp. $\Qbb$, $\Qww$, $\Qbw$) depending on the 
colors (black or white) of the two axial
vertices. We write $\gbb (\yb,\yw)$, $\gww(\yb,\yw)$, 
$\gbw(\yb,\yw)$ (resp. $\qbb (\zb,\zw)$, $\qww (\zb,\zw)$, $\qbw (\zb,\zw)$) 
for the associated series.  
\item
The families $\Gfvt$ (resp. $\Qfvt$) of 2-rooted simple
(resp. irreducible) quadrangulations 
of type \emph{face-vertex} and such that the 2-root is not incident to
the axial face. Hence 
the two roots are incident to two different faces that form an orbite
of faces of size 2 for 
the symmetry induced by the 2-root. We write $\gfvt (y)$ (resp. $\qfvt
(z)$) for their respective 
generating functions. If we take the
bicoloration into account, $\Gfvt$ (resp. $\Qfvt$) is partioned in two families
$\Gfbt$, $\Gfwt$ (resp. $\Qfbt$, $\Qfwt$) depending on the color of the axial
vertex. We write $\gfbt (\yb,\yw)$, $\gfwt(\yb,\yw)$ 
(resp. $\qfbt (\zb,\zw)$, $\qfwt (\zb,\zw)$) for the associated series.  
\item
The families $\Gfvtt$ (resp. $\Qfvtt$) of 2-rooted simple
(resp. irreducible) quadrangulations 
of type \emph{face-vertex} and such that the 2-root is incident to the
axial face. 
We write $\gfvtt (y)$ (resp. $\qfvtt (z)$) for their respective generating
functions. Similarly as above, taking the bicoloration into account,  
$\Gfvtt$ (resp. $\Qfvtt$) is partitioned in
the families  $\Gfbtt$, $\Gfwtt$ (resp. $\Qfbtt$, $\Qfwtt$) depending on the 
color of the axial
vertex. We write $\gfbtt (\yb,\yw)$, $\gfwtt(\yb,\yw)$, 
(resp. $\qfbtt (\zb,\zw)$, $\qfwtt (\zb,\zw)$) for the associated series.
\item
The families $\Gfft$ and $\Qfft$ respectively of 2-rooted simple
and irreducible quadrangulations of 
type \emph{face-face} and such that the 2-root is not incident to an
axial face. In one variable (resp. two variables), we write $\gfft (y)$ 
and $\qfft (z)$ (resp. $\gfft (\yb,\yw)$ 
and $\qfft (\zb,\zw)$) for their generating functions.
\item
The families $\Gfftt$ (resp. $\Qfftt$) of 2-rooted simple
(resp. irreducible) quadrangulations of 
type \emph{face-face} and such that the 2-root is incident to an
axial face. In one variable (resp. two variables), we write $\gfftt (y)$ 
and $\qfftt (z)$ (resp. $\gfftt (\yb,\yw)$ 
and $\qfftt (\zb,\zw)$) for their generating functions.
\end{itemize}

\begin{proposition}
The generating functions of the 5 families of 2-rooted simple
quadrangulations in one variable (resp. two variables) are
$\eta$-rational (resp. $(\eta_1,\eta_2)$-rational).
\end{proposition}
\begin{proof}
The families of 2-rooted simple quadrangulations of type vertex-vertex have 
already been 
proved in Section~\ref{sec:count2connected} to be $\eta$-rational in one 
variable and 
$(\eta_1,\eta_2)$-rational 
in two variables. 

We take the example of $\gfftt(y)$. Let $G_{f\!f,n}$
(resp. $G_{\mathfrak{f}\mathfrak{f}',n}^{(2)}$) be the number of 
2-rooted simple quadrangulations of type \emph{face-face} (resp. of
objects of $\Gfftt$) with $2n+2$ faces. 
An object of $\mathcal{G}_{\mathfrak{f}\mathfrak{f}',n}^{(2)}$ has $(2n+2)$
faces, hence it has $(8n+8)$ half-edges, 
hence it has $(4n+4)$ 2-orbites of half-edges. From this object, we
can construct a bi-2-rooted object by marking differently one of the
$(4n+4)$ orbites of half-edges. 
We obtain thus bijectively all bi-2-rooted simple quadrangulations of
type \emph{face-face} such 
that the first 2-root is incident to an axial face. Hence the family
of these objects
has cardinality $(4n+4)G_{\mathfrak{f}\mathfrak{f}',n}^{(2)}$. There is a 
second 
way to construct such a bi-2-rooted 
object, by taking a 2-rooted quadrangulation of type \emph{face-face}
and marking differently 
one of the 4 orbites of two half-edges that are incident to one of
the two  axial
faces. 

The two equivalent constructions yield the formula: 
$(4n+4)G_{\mathfrak{f}\mathfrak{f}',n}^{(2)}=4G_{f\!f,n}$. 
Hence, $\gfftt (y)=1/y\int G_{f\!f}(y)dy$. We
know that $G_{f\!f}(y)=\frac{1}{(1-3\eta(y))(1-\eta(y))}$. In addition,
$y=\eta(1-\eta)^2$, so that $dy=(1-3\eta)(1-\eta)d\eta$. 
Hence,
$\gfftt (y)=\frac{1}{\eta(1-\eta)^2}\int\frac{1}{(1-3\eta)(1-\eta)}(1-3\eta)(1-\eta)d\eta=\frac{1}{\eta(1-\eta)^2}\int
d\eta=\frac{1}{(1-\eta)^2}$. We also obtain trivially
$\gfft(y)=G_{f\!f}(y)-\gfftt(y)=\frac{2\eta}{(1-3\eta)(1-\eta)^2}$. 
The case of
$\gfvtt(y)$ can be treated similarly, giving
$\gfvtt(y)=2\eta/(1-\eta)$ and
$\gfvt(y)=G_{vf}(y)-\gfvtt(y)=\frac{4\eta}{(1-3\eta)(1-\eta)}$. 

Let us now deal with the case of two variables. We would like to obtain an 
$(\eta_1,\eta_2)$-rational expression for 
$G_{\mathfrak{f}\mathfrak{f}'}(\yb,\yw)$ from the $(\eta_1,\eta_2)$-rational 
expression of $G_{f\!f}(\yb,\yw)$. The two equivalent constructions of
bi-2-rooted objects yields the equation 
$(i+j+1)G_{\mathfrak{f}\mathfrak{f}',i,j}^{(2)}=G_{f\!f,i,j}$. 
Hence, we have to ``integrate'' $G_{f\!f}(\yb,\yw)$. Unfortunately,
unlike for the case in one variable,
there is no systematic method of integration in two variables. However,
$\gfftt$ is the only solution $f$ of the equation

\begin{equation}
\label{eq:diffeq}
\yb\frac{\partial
f}{\partial \yb}+\yw\frac{\partial f}{\partial \yw}+f=G_{\!f\!f}
\end{equation}

Hence, we just
have to guess an $(\eta_1,\eta_2)$-rational expression
$R(\eta_1,\eta_2)$  for which we have good hints that it is equal to
$\gfftt$. Then we just have to check that the $(\eta_1,\eta_2)$ rational 
expression
corresponding to $\yb\frac{\partial
R}{\partial \yb}+\yw\frac{\partial R}{\partial \yw}+R$  is equal to the 
$(\eta_1,\eta_2)$-rational
expression of $G_{\!f\!f}$. What are the hints that we have to guess the
solution ? For example, we can use the $\eta$-rational expression
$\tilde{R}(\eta)=\frac{1}{(1-\eta)^2}$ of
$\gfftt (y)$. As $\gfftt (y,y)=\gfftt(y)$ (because of the conventions given in Section~\ref{sec:conv}), 
$\eta_1(y,y)=\eta(y)$ and
$\eta_2(y,y)=\eta(y)$, a candidate $R(\eta_1,\eta_2)$ has 
to verify $R(\eta,\eta)=\tilde{R}(\eta)$. In addition $\gfftt
(\yb,\yw)$ is symmetrical in $\yb$ and $\yw$. As
$\eta_1(\yb,\yw)=\eta_2(\yw,\yb)$, a candidate has to be
symmetric in $\eta_1$ and $\eta_2$ (because such a candidate is symmetrical 
in $\yb$ and $\yw$). These two hints lead us to guess that
$R(\eta_1,\eta_2)=\frac{1}{(1-\eta_1)(1-\eta_2)}$. It turns out that
this candidate verifies Equation~(\ref{eq:diffeq}), hence $\gfftt =
\frac{1}{(1-\eta_1)(1-\eta_2)}$. We also get an
$(\eta_1,\eta_2)$-rational expression for $\gfft$, using
$\gfft=G_{f\!f}-\gfftt$. Similarly, we can also guess and check an
$(\eta_1,\eta_2)$-rational expression for $\gfbtt$. We find
$\gfbtt=\frac{\eta_1}{1-\eta_2}$, and find also an
$(\eta_1,\eta_2)$-rational expression for $\gfbt$ using the relation 
$\gfbt=G_{bf}-\gfbtt$. We also find $(\eta_1,\eta_2)$-rational
expressions for $\gfwtt$ and $\gfwt$, observing that
$\gfwtt(\yb,\yw)=\gfbtt(\yw,\yb)$ and $\gfwt(\yb,\yw)=\gfbt(\yw,\yb)$.
Hence, one obtains $(\eta_1,\eta_2)$-rational expressions of $\gfbtt$ and $\gfbt$  by substituting $(\eta_1,\eta_2)$ by
$(\eta_2,\eta_1)$ in the respective $(\eta_1,\eta_2)$-rational
expressions of $\gfbtt$ and $\gfbt$.
\end{proof}

In order to find expressions of the 5 generating functions of 2-rooted 
irreducible quadrangulations, we will relate them to the 5 generating 
functions of 2-rooted simple quadrangulations (which we have obtained) by 
a triangular sytem of 5 equations: for each family of 2-rooted simple 
quadrangulations, an equation is derived from a canonical decomposition of 
the objects of the family, with an irreducible quadrangulation at 
the ``core'' of the decomposition. 

To derive these equations, we will also need the two following auxiliary 
families of 2-rooted simple
quadrangulations, 
as intermediates of calculation:
\begin{itemize}
\item
Let $\mathcal{L}$ be the family of 2-rooted simple quadrangulations of
type \emph{face-face} 
such that the 2-root is incident to an axial face and such that the
root node of the induced 
rooted quadrangulation (the induced rooted quadrangulation is
obtained by keeping 
only the mark of one of the two roots) is not a vertical
axis-node. We write $L(y)$ (resp. $L(\yb,\yw)$) for the 
series counting $\mathcal{L}$ in one variable (resp. two variables). In two variables, we denote by $^tL$ the series $^tL(\yb,\yw):=L(\yw,\yb)$.
\item
Let $\mathcal{K}$ be the family of 2-rooted simple quadrangulations of
type \emph{face-vertex} such that the 2-root is incident to the
axial face and 
such that the root node of the induced rooted quadrangulation is not a
vertical axis-node. 
We write $K(y)$ for the series counting $\mathcal{K}$ in one variable. In
two variables, we write $K_{b}(\yb,\yw)$
(resp. $K_{w}(\yb,\yw)$) for the series counting objects of
$\mathcal{K}$ whose axial vertex is black (resp. white). We also denote by $^tK_b$ and $^tK_w$ the series $^tK_b(\yb,\yw):=K_b(\yw,\yb)$ and $^tK_w(\yb,\yw)=K_w(\yw,\yb)$.
\end{itemize}

\subsubsection{Calculation of the generating functions of $\mathcal{K}$ and 
$\mathcal{L}$.}
\label{sec:h}
\begin{figure}
\begin{center}
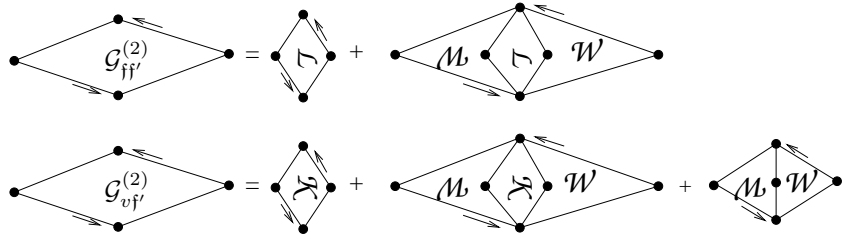
\caption{The two decompositions involving respectively the objects of 
$\mathcal{L}$ and $\mathcal{K}$.}
\label{figure:declh}
\end{center}
\end{figure}

There is an easy decomposition  of the objects of $\Gfftt$
(resp. $\Gfvtt$) using the 
objects of $\mathcal{L}$ (resp. $\mathcal{K}$). These decompositions
are performed by 
looking if the root node of the decomposition-tree is an horizontal
axis-node or not, 
see Figure~\ref{figure:declh}. These two decompositions give the 
following equations, ensuring that the series of $\mathcal{L}$ and 
$\mathcal{K}$ in one variable (resp. two variables) are
$\eta$-rational (resp. $(\eta_1,\eta_2)$-rational).
\begin{equation} 
\gfftt (y)=L(y)(\ovw+1)
\end{equation}
\begin{equation}
\label{equation:l2}
\gfftt \left(\yb,\yw\right)=L\frac{W}{\yb}+L
\end{equation}

\begin{equation}
\label{equation:h}
\gfvtt (y)=K(y)(\ovw +1)+\ovw
\end{equation}
\begin{equation}
\left\{
\begin{array}{rcl}
\gfbtt (\yb,\yw)&=&K_{b}\frac{W}{\yb}+K_{b} \\
\gfwtt (\yb,\yw)&=&K_{w}\frac{W}{\yb}+K_{w}+\frac{W}{\yb}
\end{array}
\right.
\end{equation}

\subsubsection{Calculation of the generating function of $\Qfftt$.}
\label{section:q44tt}
An object of $\mathcal{L}$ has a simple decomposition: it is either
the trivial quadrangulation (i.e. with only two faces), or the root
node of its decomposition tree
is a vertical axis-node or is
an irreducible 
quadrangulation.

This decomposition yields in one variable (resp two variables) the
following equations:
\begin{equation}
\label{equation:q44tt}
L(y)=1+L(y)\ovw +\qfftt \left( \ovw\right)\gfftt (y)
\end{equation}
\begin{equation}
L(\yb,\yw)=1+^tL\frac{W}{\yw}+\qfftt\yW\gfftt(\yb,\yw).
\end{equation}

These equations ensure that $\qfftt \left( \ovw\right)$ and
$\qfftt\yW$ are respectively $\eta$-rational and $(\eta_1,\eta_2)$-rational.

\subsubsection{Calculation of the generating function of $\Qfvtt$.}
\label{section:q42tt}
Similarly as above, we perform here a decomposition of an object of
$\mathcal{K}$ by 
looking if the root-node of its decomposition-tree is a vertical
axis-node or an irreducible 
quadrangulation. In two variables, we perform the same decomposition,
but we distinguish whether the axial vertex is black or white, see
Figure~\ref{figure:dech} for the decomposition of an object of
$\mathcal{K}_{b}$. This yields the equations:

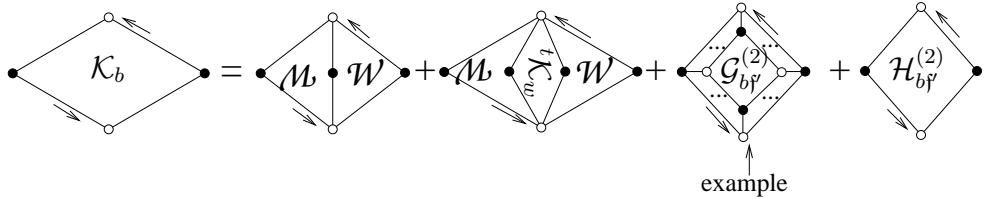
\begin{figure}
\begin{center}
\input{Figures/dech.pstex_t}
\caption{Decomposition of an objet of $\mathcal{K}_{b}$.}
\label{figure:dech}
\end{center}
\end{figure}

\begin{equation}
\label{equation:q42tt}
K(y)=\ovw+K(y)\ovw+\qfftt \left( \ovw\right)\gfvtt
(y)+\qfvtt \left( \ovw\right)
\end{equation}
\begin{equation}
\label{eq:q42tt2}
\left\{
\begin{array}{rcl}
K_{b}(\yb,\yw)&=&\frac{W}{\yw}+^tK_{w}\frac{W}{\yw}+\qfftt\yW\gfbtt (\yb,\yw)+\qfbtt\yW\\
K_{w}(\yb,\yw)&=&^tK_{b}\frac{W}{\yw}+\qfftt\yW\gfwtt
(\yb,\yw)+\qfwtt\yW
\end{array}
\right.
\end{equation}
We already know from Sections~\ref{sec:h} and~\ref{section:q44tt} that 
$\qfftt \left( \ovw\right)$ and $K(y)$ are $\eta$-rational and that $\qfftt\yW$, $K_{b}(\yb,\yw)$ 
and
$K_{w}(\yb,\yw)$ are $(\eta_1,\eta_2)$-rational. Hence
Equation~(\ref{equation:q42tt}) yields an $\eta$-rational
expression for $\qfvtt \left( \ovw\right)$ and
Equation~(\ref{eq:q42tt2}) yields 
$(\eta_1,\eta_2)$-rational expressions for  $\qfbtt\yW$ and
$\qfwtt\yW$.

\subsubsection{The trunk of the decomposition-tree and the core-node.}
Before proceeding further, we need a better understanding of how the 
symmetry of a 2-rooted quadrangulation is 
repercuted  on its decomposition-tree. 

Let $Q$ be a 2-rooted quadrangulation whose decomposition-tree has at 
least one node fixed by the induced symmetry, and whose 2-root is not 
incident to an axial face. As the axis of the rotation-symmetry passes 
by all nodes fixed by the symmetry, these nodes form a chain 
$x_0,\ldots,x_k$, with $k\geq 0$, of 2-rooted
quadrangulations, $x_i$ and $x_{i-1}$ being connected 
at a common axial face.  The chain $x_0,\ldots,x_k$ 
is called the \emph{trunk} of the 
decomposition tree. Observe that $x_i$ has type face-face 
for $1\leq i\leq k-1$, and that the two axial cells of $Q$ 
are the axial cell of $x_0$ and the axial face of $x_k$ not 
involved in an interconnection.
 
As the 2-root is not incident to an axial face of $Q$, there exists a 
node $x_Q$ of the trunk such that two identical quadrangulations, 
carrying each a root of $Q$, are connected at each face of a 2-orbit of 
faces of $x_Q$. The node $x_Q$ is called the \emph{core-node} of $Q$.

\begin{proposition}
\label{prop:center_tree}
A 2-rooted
simple quadrangulation $Q$ with at least one axial face has at least one node 
of its decomposition-tree fixed by the induced symmetry.
\end{proposition}
\begin{proof}
Let $Q$ be a 2-rooted quadrangulation whose decomposition tree has no node 
fixed by the symmetry induced by the 2-root.
Hence the centre of the decomposition tree of $Q$  has to be an edge, and 
the symmetry has to exchange the two nodes incident to the edge. 
Geometrically, an edge of the decomposition tree corresponds
to an interconnection-face $f$ connecting two nodes of the tree at two
faces of these nodes. The symmetry can not exchange two vertices of $Q$ with 
different colors when $Q$ is bicolored. Hence, 
the symmetry necessarily turns over $f$ around two vertices $s_1$ and
$s_2$ that are diagonally 
opposed in $f$ and that are hence the two axial cells of $Q$. Thus $Q$
is necessarily of 
type \emph{vertex-vertex}.
\end{proof} 

\subsubsection{Calculation of the generating function of $\Qfft$.}
\label{section:q44t}
The idea is to perform a decomposition of an object of $\Gfft$
according to several cases. 
This decomposition yields an equation linking the generating functions of 
the families $\Gfftt$, $\Gfft$ and $\Qfft$.

Let $Q\in\Gfft$. The rotation-symmetry of order 2 of
$Q$ induces a 
symmetry of its decomposition-tree. According to 
Proposition~\ref{prop:center_tree}, $Q$ has at least one node of 
its decomposition-tree fixed by the symmetry. Hence the definitions 
of trunk and of core-node apply for $Q$. Then there are two possibilities 
for the core-node $x_Q$ of $Q$: either it is an axis-node or it is an 
irreducible quadrangulation.

\textbf{The core-node is an axis-node.} 
We present here a construction of composed objects that correspond 
bijectively to objects
of $\Gfft$ whose repercussion-node is an axis-node with
$2d+2$ faces. Take the
2-rooted axis-node $A$ of type face-face, with 4 faces, and such that
its 2-root is incident to an axial face with each of the two roots pointing
toward each of the two extremal vertices of $A$.

Then take two copies of an object of $\mathcal{W}'$ whose vertical
multiplicity at the root node of its decomposition tree is $d$. Root each of 
these 
two copies on each of the two non rooted faces of $A$, superposing the edge
carrying the primary 
root of each copy of $\mathcal{W}'$ with the edge
carrying each of the two roots of $A$. Then, on each axial
face of $A$, 
root an object of $\mathcal{L}$ in a canonical way and ``vertically''.
%\begin{figure}
%\begin{center}
%\includegraphics{Figures/decaxismil.eps}
%\caption{Construction of an object of $\Gfft$ such that its repercussion-node is an axis-node.}
%\label{figure:decaxismil}
%\end{center}
%\end{figure}
Finally keep only the mark of the secondary roots of each of the two
copies of the 
objects of $\mathcal{W}'$. 

This construction allows us to obtain exactly twice all object of $\Gfftt$ 
whose core-node is an axis-node (the inverse
construction consists in guessing which axial face of $A$ was carrying
the 2-root, which gives two choices).

Let $r_{\mathfrak{f}\mathfrak{f}}^{(2)}(y)$ be the generating function of
objects of $\Gfftt$ 
whose core-node is an axis-node. From the construction explained above, we have

\begin{equation}
r_{\mathfrak{f}\mathfrak{f}}^{(2)}(y)=2\frac{C(y)}{y}L(y)^2
\label{equation:q44taxis}
\end{equation}

In two variables, we have: 
\begin{equation}
r_{\mathfrak{f}\mathfrak{f}}^{(2)}(\yb,\yw)=\frac{C}{\yw} {^tL^2}+\frac{C}{\yb}L^2
\end{equation}

\textbf{The core-node is an irreducible quadrangulation.}
The method 
is still to construct objects of $\Gfftt$ whose core-node is an irreducible 
quadrangulation as composed objects. 
Take $Q\in\Qfftt$ and take two copies of an object of
$\mathcal{W}'$. Root canonically each of these copies on each of the two rooted
faces of $Q$. 
Then, for each 2-orbite of non rooted faces of $Q$, take two copies
of an object 
of $\mathcal{W}$ and root canonically each copy on each face of the
orbite. 
Then root canonically an object of $\Gfftt$ on each axial face of
$Q$. 
Finally, keep only the mark of the secondary roots of the two copies of the
object 
of $\mathcal{W}'$. Similarly as in Section~\ref{sec:quadsimpb},
each object of $\Gfft$ whose core-node is an irreducible
quadrangulation is obtained exactly 4 times by this construction. 
We obtain the following series for these objects, in one and in two variables:

$$
\frac{C(y)}{W(y)}\qfft (\ovw)\gfftt (y)^2
$$
$$
\frac{C}{W}\qfft\yW \gfftt
(\yb,\yw)^2
$$
\textbf{Obtaining the equations.} Summing the two contributions for the
decomposition of an object of $\Gfft$, we obtain the following equations , in
one and two variables:

\begin{equation}
\gfft (y)=2\frac{C(y)}{y}L(y)^2+\frac{C(y)}{W(y)}\qfft (\ovw)\gfftt (y)^2
\label{equation:q44t}
\end{equation}

\begin{equation}
\label{eq:qfacevariable}
G_{\mathfrak{f}\mathfrak{f}}^{(2)}(\yb,\yw)=\frac{C}{\yw} {^tL^2}+\frac{C}{\yb}L^2+\frac{C}{W}H_{\mathfrak{f}\mathfrak{f}}^{(2)}\yW \gfftt (\yb,\yw)^2.
\end{equation}

We deduce from Equation~(\ref{equation:q44t}) that $\qfft \left(
\ovw\right)$ is $\eta$-rational, because all other series
appearing in Equation~(\ref{equation:q44t}) are already known and
$\eta$-rational. Similarly, using Equation~(\ref{eq:qfacevariable}), we
can obtain an $(\eta_1,\eta_2)$-rational expression for $\qfft \yW$.

\subsubsection{Calculation of the generating function of $\Qfvt$.}
\label{section:q42t}
In this section, we find an equation linking the families $\Gfvt$, $\Gfvtt$,
$\Gfftt$, $\Qfft$ and $\Qfvt$. 
This equation is derived from a construction of objects of $\Gfvt$
as composed objects, distinguishing several cases.

Let $Q\in\Gfvt$. According to 
Proposition~\ref{prop:center_tree}, the decomposition-tree of $Q$ has at 
least one node fixed by the symmetry induced by the 2-root, so that the 
definitions of trunk and core-node apply. Let
$x_0,\ldots ,x_k$ be the trunk of $Q$. As $Q$ has type face-vertex, one 
extremity of the
trunk, say $x_0$, has type 
\emph{face-vertex}. As in Section~\ref{section:q44t}, we have two 
possibilities for the core-node $x_Q$: either $x_Q$ is an axis-node or $x_Q$ 
is an 
irreducible quadrangulation. 

\textbf{The core-node is an axis-node.} 
Two subcases can arise:
\begin{itemize}
\item
The axis-node $x_Q$ is not the node $x_0$ of the trunk, i.e. $x_Q$ has 
type \emph{face-face}.
\item
The axis-node $x_Q$ is the node $x_0$ of the trunk, i.e. $x_Q$ has type 
\emph{face-vertex}.
\end{itemize}

We treat the first case with a construction similar to the one that has 
lead to
Equation~(\ref{equation:q44taxis}) in Section~\ref{section:q44t}. The
only difference is that, at the end, we do not root two objects of
$\mathcal{L}$ on 
the axial faces of the axis-node $A$, but we root ``vertically'' an
object of $\mathcal{K}$ 
on the rooted axial face and an object of $\mathcal{L}$ on the second
axial face of $A$. With this construction we obtain bijectively all
obects of $\Qfvt$ such that $x_{Q}$ is an axis-node of type face-face.   
A similar construction can be performed by
taking the bicoloration of vertices (and in particular the color of
the axial vertex of the object of $\mathcal{K}$) into account. This yields
the following generating functions for this case,
respectively in one and two variables
 (the generating functions in two variables are those of objects whose axial 
vertex is black, as the
generating function
for objects whose axial vertex is white can be deduced by exchanging the 
variables):

$$
4\frac{C(y)}{y}K(y)L(y)
$$
$$
2^tL^tK_{w}\frac{C}{\yw}+2LK_{b}\frac{C}{\yb}
$$

The second subcase is treated by performing a quite analogous construction
of composed objects. Take the 2-rooted axis-node $A$ with 3 faces and
one axial face, such that its two roots are incident to the axial face and 
point toward each of the two extremal vertices of $A$. Then root two copies of an 
object of
$\mathcal{W}'$ on each of the two non axial faces of $A$, superposing
the root edge of each of the two copies with each of the two root
edges of $A$. This construction gives rise bijectively to all objects of
$\Qfvt$ whose core-node is an axis-node with type face-vertex.  
We obtain the following generating functions respectively in one and in two
variables for these objects (once again, in two variables, we only
consider objects whose axial vertex is black):
$$
4\frac{C(y)}{y}L(y)
$$
$$
2^tL\frac{C}{\yw}
$$

%Finally, summing the contributions of the two subcases, we obtain
%respectively the generating function $r_{v\mathfrak{f}}(y)$ and $r_{b\mathfrak{f}}(\yb,\yw)$ of the objects of $\Gfvt$ and of
%$\Gfbt$ whose core-node is an axis-node:

%$$
%r_{v\mathfrak{f}}(y)=4\frac{C(y)}{y}l(y)\left( 1+h(y)\right)
%\label{equation:q42taxis}
%$$
%$$
%r_{b\mathfrak{f}}(\yb,\yw)=2^t\!l^t\!h_{\circ}\frac{C}{\yw}+2lh_{\bullet}\frac{C}{\yb}+2^t\!l\frac{C}{\yw}
%$$

\textbf{The core-node is an irreducible quadrangulation.} Two subcases can 
also arise here  depending on $x_Q$
having type \emph{face-face} or type \emph{face-vertex}. We perform a
similar composition construction as in Section~\ref{section:q44t} for objects 
of $\Gfft$ whose core-node is irreducible. The only
difference is that, at the end, we do not root two objects of $\Gfftt$
on 
the axial faces of the object of $\Qfft$, but we root an object of
$\Gfftt$ 
on one axial face and an object of $\Gfvtt$ on the other axial face of
the irreducible quadrangulation. 
We obtain for this case the following generating functions respectively in one
and two variables:
$$
2\frac{C(y)}{W(y)}\qfft \left(\ovw\right)\gfftt (y)\gfvtt(y)
\label{equation:q42tirreductiblefaceface}
$$
$$
2\frac{C}{W}\qfft \yW \gfftt (\yb,\yw)\gfbtt (\yb,\yw)
$$

Then we treat the case where $x_Q$ is an irreducible quadrangulation
and has type \emph{face-vertex}. 
This case is treated similarly as the preceding case, with the difference 
that  we
take an irreducible 
quadrangulation $Q$ with type \emph{face-vertex} at the beginning of
the construction. 
In addition, at the end, we root an object of $\Gfftt$ on the
unique axial face of $Q$. 
We obtain for this case the generating functions
$$
\frac{C(y)}{W(y)}\qfvt \left(\ovw\right)\gfftt (y)
\label{equation:q42tirreductiblefacevertex}
$$
$$
\frac{C}{W}\qfbt \yW \gfftt (\yb,\yw).
$$

Finally, we obtain the following equations
corresponding to all 
cases for the decomposition of an object of $\Gfvt$ (resp. of $\Gfbt$):
\begin{eqnarray}
\gfvt (z)&=&4\frac{C(y)}{y}\left( 1+K(y)\right)+\frac{C(y)}{W(y)}\qfvt\left(\ovw\right)\gfftt (y)\label{equation:q42t}\\
&&+2\frac{C(y)}{W(y)}\qfft (\ovw)\gfftt (y)\gfvtt (y)\nonumber
\end{eqnarray}

\begin{eqnarray}
\gfbt (\yb,\yw)&=&2^tL^tK_{w}\frac{C}{\yw}+2L\ \!K_{b}\frac{C}{\yb}+2^tL\frac{C}{\yw}+\frac{C}{W}\qfbt \yW \gfftt (\yb,\yw)\nonumber\\
&&+2\frac{C}{W}\qfft \yW
\gfftt (\yb,\yw)\gfbtt (\yb,\yw)\label{eq:qfbt}
\end{eqnarray}

The only unknown generating function in
Equation~(\ref{equation:q42t}) is $\qfvt (W(y)/y)$, as all other series
appearing in this equation are already known to be
$\eta$-rational. Hence we obtain from this equation an $\eta$-rational
expression for $\qfvt (W(y)/y)$.
Similarly, we obtain from Equation~(\ref{eq:qfbt}) an
$(\eta_1,\eta_2)$-rational expression for $\qfbt \yW$. Finally, we
also obtain easily an $(\eta_1,\eta_2)$-rational expression for
$\qfwt \yW$. Indeed, we have clearly
$\qfwt(\zb,\zw)=\qfbt(\zw,\zb)$. Observing that
$\eta_1(\yb,\yw)=\eta_2(\yw,\yb)$ and
$\eta_2(\yb,\yw)=\eta_1(\yw,\yb)$, we obtain an
$(\eta_1,\eta_2)$-rational expression of $\qfwt \yW$ by exchanging $\eta_1$ 
and $\eta_2$ in
the $(\eta_1,\eta_2)$-rational expression of $\qfbt \yW$.

\subsubsection{Calculation of the generating function of $\Gvv$.}
\label{section:q22}

It now remains to find the generating function of the
family $\Qvv$. 
To do this, we construct objects of $\Gvv$ as composed objects  by
distinguishing several cases. 
We will obtain thus an equation linking $\qvv (z)$ and the generating
functions of the 
other families of 2-rooted simple and irreducible quadrangulations.

\textbf{The different cases for the repercussion of the symmetry.}
Let $Q$ be an object of $\Gvv$. This time a lot of cases can arise for
the repercussion of the 
symmetry induced by the 2-root of $Q$ on the
decomposition-tree of $Q$. Indeed, it can happen that no node of the 
decomposition tree is fixed by the symmetry. As we have seen in the proof 
of Proposition~\ref{prop:centre}, this corresponds to the situation where 
the centre of the 
decomposition tree is an edge whose two extremities are exchanged by the 
symmetry: on the quadrangulation, this edge corresponds to an 
interconnection-face that is turned over by the symmetry. This situation 
is treated in Case 2.

Otherwise we consider the \emph{trunk} of $Q$  (i.e. the chain of nodes
$x_0,\ldots ,x_k$ that are 
fixed by the symmetry) and the core-node $x_Q$. Two cases can arise for 
the trunk: $k=0$
i.e. the trunk has 
only one node $x_0$ of type \emph{vertex-vertex} and then $x_Q=x_0$;
or $k\geq 1$, which 
implies that $x_0$ and $x_k$ have type \emph{face-vertex} and that $x_i$ 
has type \emph{face-face} $\forall 1\leq i\leq k-1$. Thus, if $k\geq 1$, 
$x_Q$ can have type
\emph{face-vertex} if it is 
an extremity of the trunk, otherwise it has type \emph{face-face}. In
addition, we have each 
time to distinguish whether $x_Q$ is an axis-node or an irreducible
quadrangulation. Moreover 
we have to take care of the fact that, if $x_Q$ is an axis-node, there
are two distinct ways 
for $x_Q$ to have type \emph{vertex-vertex} (treated respectively in cases 3 
and 4).

We distinguish the following cases:
\begin{enumerate}
\item
The quadrangulation $Q$ has only two faces and its 2-root induces a
symmetry of $Q$ 
with two axial vertices. This case can be considered as degenerated
by saying that the 
decomposition-tree of $Q$ consists only of a degenerated axis-node with 2 
faces.
\item
The centre of the decomposition-tree of $Q$ is an edge, and this
edge is turned 
over by the symmetry. This edge corresponds to an interconnection-face
$f$ connecting 
two nodes of the tree. Hence the two axial vertices of $Q$ are two
diagonally opposed 
vertices of $f$. Moreover it is easy to see that the two nodes connected 
at $f$ are irreducible quadrangulations. Indeed, assume that these nodes 
$n_1$ and $n_2$
are 
axis-nodes. As we have seen in Section~\ref{section:treedecompsecond}, two 
incident axis-nodes are 
stretched in perpendicular
directions: 
i.e. if $n_1$ is ``horizontal'' then $n_2$ is ``vertical''. The
symmetry turns over 
the tree around the face $f$ and sends $n_1$ to the place where $n_2$
was. But $n_1$ 
clearly remains ``horizontal'' after this turn-over, and takes the
place of
$n_2$, which is vertical. Hence this turn over operation can not let
$Q$ invariant, so that we have a contradiction.
Hence the two connected 
nodes are irreducible quadrangulations.
\item
The core-node $x_Q$ is an axis-node and the rotation-axis of the symmetry
of $Q$ induced by 
the 2-root intersects $x_Q$ at two diagonally opposed vertices in the 
equatorial plane of $x_Q$.
\item
The core-node $x_Q$ is an axis-node and the rotation-axis of the symmetry
of $Q$ induced by 
the 2-root intersects $x_Q$ at its   two extremal vertices.
\item
The core-node $x_Q$ is an axis-node and the rotation-axis of the symmetry
of $Q$ induced by 
the 2-root intersects $x_Q$ at  a vertex $v$ of the equatorial plane of $x_Q$ 
and
at 
the centre of the face diametrically opposed to $v$
in the equatorial 
plane of $x_Q$. This corresponds to the case where $x_Q$ is an axis-node of 
type \emph{face-vertex}.
\item
The core-node $x_Q$ is an axis-node and the rotation-axis of the symmetry
of $Q$ 
induced by the 2-root  intersects $x_Q$ at the centres of two diametrically
opposed faces of $x_Q$. This corresponds to the case where $x_Q$ is an 
axis-node of type \emph{face-face}.
\item
The core-node $x_Q$ is an irreducible quadrangulation that has type 
\emph{face-face} for the symmetry induced by the 2-root.
\item
The core-node $x_Q$ is an irreducible quadrangulation that has type 
\emph{face-vertex} for the symmetry induced by the 2-root.
\item
The core-node $x_Q$ is an irreducible quadrangulation that has type 
\emph{vertex-vertex} for the symmetry induced by the 2-root.
\end{enumerate}

\textbf{The core-node is a (possibly degenerated) axis-node.}
Cases 1, 2 and 3 can be treated together by constructing the
following 
composed objects: let $A$ be the 2-rooted simple quadrangulation
with two 
faces and type \emph{vertex-vertex} and such that the two roots point
towards 
the same vertex of $A$. We call this vertex the root-vertex of $A$. We
can 
consider $A$ as a sort of degenerated
axis-node 
with only two faces. Hence we can perform the same construction as in 
Section~\ref{sec:quadsimp_a}. Here this construction comes
down to 
taking two copies of an object $Q'$ of $\widehat{\mathcal{J}'}$ and rooting
each copy on each of the two faces of $A$ so that the primary roots of the
two copies 
point towards the root vertex of $A$. By definition of $\mathcal{J}$, $Q'$ is
not stretched 
horizontally at its root face (the face incident to the primary root). 
If $Q'$ is not vertically stretched, then $Q'$ is either a
trivial 
quadrangulation with two faces, corresponding to Case 1, or
$Q'$ 
is an irreducible quadrangulation, corresponding to Case
2. 
The situation where $Q'$ is vertically streched at its root face corresponds 
to Case 3.

Finally, the sum of the contributions of Case 1, Case 2 and Case 3 is
\[
2\frac{B(y)}{y}.
\]

Case 4 can also be treated similarly as in
Section~\ref{sec:quadsimp_a}. 
However we have to take care of the fact that the axis-node at the core of 
the decomposition-tree 
is a real axis-node, so that it has at least 3 faces. Hence it has
$2\cdot l$ faces with $l\geq 2$. 
Thus Case 4 gives the generating function
\[
2\frac{B(y)}{y}\frac{\ovg}{1-\ovg}
\]

Finally, the sum of the contributions of Case 1, Case 2, Case 3 and Case 4 is

$$
2\frac{B(y)}{y}\frac{1}{1-J(y)/y}.
\label{equation:q22axisdegenerated}
$$

Case 5, where the core-node $x_Q$ is an axis-node of type \emph{face-vertex},
can be treated similarly as in Section~\ref{section:q42t}. For this
construction, we took a 2-rooted axis-node $A$ with 3 faces and with type
\emph{face-vertex} at the beginning of the construction. The only
difference is that, at the end of the construction, we do not root
``vertically'' 
an object of $\mathcal{L}$ on the axial face of $A$ but an object of
$\mathcal{K}$. 
We obtain for Case 5 the generating function

$$
4\frac{C(y)}{y}K(y).
\label{equation:q22axisfacevertex}
$$

Case 6 can also be treated similarly as in the construction of
Section~\ref{section:q44t} 
where we took a 2-rooted axis-node $A$ with 4 faces and 2 axial
faces. The only difference is 
that, at the end of the construction, we root two objects of
$\mathcal{K}$ ``vertically'' on 
the two axial faces of $A$. Hence Case 6 gives the generating function

$$
2\frac{C(y)}{y}K(y)^2.
\label{equation:q22axisfaceface}
$$

Finally we can group the 6 first cases in a generating function
$r_{vv}^{(2)}(y)$ of objects of $\Gvv$ 
whose core-node is a (possibly degenerated) axis-node:

$$
r_{vv}^{(2)}(y)=2\frac{B(y)}{y}\frac{1}{1-J(y)/y}+4\frac{C(y)}{y}K(y)+2\frac{C(y)}{y}K(y)^2
\label{equation:q22axis}
$$

\textbf{The core-node is an irreducible quadrangulation.}
Cases 7, 8 and 9 correspond to objects of $\Gvv$ whose core-node is an 
irreducible quadrangulation. These cases can be treated by
constructing composed objects  
from a 2-rooted irreducible quadrangulation that has type
\emph{face-face} for Case 7, 
type \emph{face-vertex} for Case 8 and type \emph{vertex-vertex} for Case 9.

Hence, the generating function for Case 7 is
$$
\frac{C(y)}{W(y)}\qfft \left( \ovw\right)\gfvtt (y)^2,
\label{equation:q22irreductiblefaceface}
$$
the generating function for Case 8 is
$$
\frac{C(y)}{W(y)}\qfvt \left( \ovw\right)\gfvtt (y),
\label{equation:q22irreductiblefacevertex}
$$
and the generating function for Case 9 is
$$
\frac{C(y)}{W(y)}\qvv \left( \ovw\right).
\label{equation:q22irreductiblevertexvertex}
$$

\textbf{Obtaining the equation.} Finally we obtain the following
equation corresponding 
to the different ways to construct objects of $\Gvv$ as composed objects:

\begin{eqnarray}
G_{vv}^{(2)}(y)&=&r_{vv}^{(2)}(y)+\frac{C(y)}{W(y)}\qfft \left( \ovw\right)\gfvtt (y)^2\label{equation:q22}\\
&&+\frac{C(y)}{W(y)}\qfvt \left( \ovw\right)\gfvtt (y)+\frac{C(y)}{W(y)}\qvv \left( \ovw\right).\nonumber
\end{eqnarray}

Except for $\qvv
\left( \ovw\right)$, all generating functions of this equation are
known and are $\eta$-rational. 
Hence we obtain from this equation an $\eta$-rational expression for $\qvv
\left( \ovw\right)$.

Similarly, for two variables, the same decomposition yields the
following equations depending on the colors of the two axial vertices:

\begin{eqnarray}
\gww\yW
&=&r_{ww}(\yb,\yw)+\frac{\yb}{\yw}\frac{C}{W}\qfft\yW
\gfwtt(\yb,\yw) ^2\nonumber \\
&&+\frac{\yb}{\yw}\frac{C}{W}\qfwt\yW
\gfwtt(\yb,\yw)+\frac{C}{W}\qww\yW,\nonumber
\end{eqnarray}
where 
$$
r_{ww}(\yb,\yw)=\frac{B}{\yw}\frac{1}{1-J/\yw}+2\frac{C}{\yw}K_{w}+\frac{C}{\yw}K_{w}^2+\frac{\yb C}{\yw ^2}{^tK_{b}^2}.
$$

\begin{eqnarray*}
\gbw\yW&=&r_{bw}(\yb,\yw)+2\frac{C}{W}\qfft \yW \gfbtt
(\yb,\yw)\gfwtt (\yb,\yw)\\
&&+\frac{C}{W}\qfbt \yW \gfwtt (\yb,\yw)\\
&&+\frac{C}{W}\qfwt \yW \gfbtt
(\yb,\yw)+\frac{C}{W}\qbw\yW,
\end{eqnarray*}
where 
$$ r_{bw}(\yb,\yw)=2\frac{C}{\yb}K_{b}+2\frac{C}{\yw}{^tK_{b}}+2\frac{C}{\yb}K_{b}K_{w}+2\frac{C}{\yw}{^tK_{b}^tK_{w}}.
$$
We conclude from these equations that $\qww \yW$ and $\qbw
\yW$ are $(\eta_1,\eta_2)$-rational. Moreover, 
$\qbb\yW=\qww \left( W/\yb,W/\yw\right)$, so that $\qbb\yW$ has also 
an $(\eta_1,\eta_2)$-rational expression, obtained by exchanging 
$\eta_1$ and $\eta_2$ in the $(\eta_1,\eta_2)$-rational expression of 
$\qww\yW$.

\subsubsection{Conclusion}
We have proved that all series of 2-rooted
quadrangulations in one variable, composed with $W(y)/y$, are 
$\eta$-rational and all
series of 2-rooted quadrangulations in two variables, composed with
$\yW$, are $(\eta_1,\eta_2)$-rational. This allows us to state the
following lemma, which completes Lemma~\ref{lemma:gammarational} (see
the proof of Lemma~\ref{lemma:gammarational} for the transition
between $(\eta_1,\eta_2)$-rational and $(\gamma_1,\gamma_2)$-rational):

\begin{lemma}
\label{lemma:2rootedirred}
All series of 2-rooted irreducible quadrangulations are $\gamma$-rational 
in one
variable and are $(\gamma_1,\gamma_2)$-rational in two variables. 
\end{lemma}

\subsection{Complexity result}
Burnside's formula for 3-connected maps can be formulated as follows: let $h_n$ be the number of unrooted 3-connected maps with $n$ edges and $h_{ij}$ be the number of unrooted 3-connected maps with $i+1$ vertices and $j+1$ edges. Then
$$
\sum_n 2nh_n z^n=H(z)+zH_{v\mathfrak{f}'}(z^2)+zH_{v\mathfrak{f}}(z^2)+z^2H_{\mathfrak{f}\mathfrak{f}'}(z^2)+z^2H_{\mathfrak{f}\mathfrak{f}}(z^2)+H_{vv}^{(2)}(z^2)+\sum_{k\geq 3}\phi(k)H_{vv}^{\geq 3}(z^k),
$$
\begin{eqnarray*}
\sum_{i,j} 2(i+j)h_{ij} \zb^i\zw^j&=&H(\zb,\zw)+\zw H_{b\mathfrak{f}'}(\zb^2,\zw^2)+\zw H_{b\mathfrak{f}}(\zb^2,\zw^2)+\zb H_{w\mathfrak{f}'}(\zb^2,\zw^2)\\
&&+\zb H_{w\mathfrak{f}}(\zb^2,\zw^2)+\zb\zw H_{\mathfrak{f}\mathfrak{f}'}(\zb^2,\zw^2)+\zb\zw H_{\mathfrak{f}\mathfrak{f}}(\zb^2,\zw^2)\\
&&+\frac{\zb}{\zw}H_{bb}^{(2)}(\zb^2,\zw^2)+\frac{\zw}{\zb}H_{ww}^{(2)}(\zb^2,\zw^2)+H_{bw}^{(2)}(\zb^2,\zw^2)\\
&&+\sum_{k\geq 3}\phi(k)\left(\frac{\zb}{\zw}H_{bb}^{\geq 3}(\zb^k,\zw^k)+\frac{\zw}{\zb}H_{ww}^{\geq 3}(\zb^k,\zw^k)+H_{bw}^{\geq 3}(\zb^k,\zw^k)\right).
\end{eqnarray*}
Lemma~\ref{lemma:gammarational} and Lemma~\ref{lemma:2rootedirred} imply
the following result of complexity for the enumeration of unrooted
3-connected maps:

\begin{lemma}
\label{lemma:complexity3}
The $N$ initial coefficients counting unrooted 3-connected maps
according to their number of edges can be computed with
$\mathcal{O}(N\log (N))$ operations.

The table of initial coefficients with indices $(i,j)$ and $i+j\leq N$
counting unrooted 3-connected maps according to their number of vertices and
faces can be computed with $\mathcal{O}(N^2)$ operations.
\end{lemma}
\begin{proof}
The proof is similar to that of Lemma~\ref{lemma:complexity2}, i.e. 
we use the property of D-finiteness of the series of rooted and $k$-rooted 3-connected maps, following from 
the fact that these series are algebraic. Indeed, in one variable 
they are in the algebraic extension of
$\gamma$, and in two variables they are in the algebraic extension of $(\gamma_1,\gamma_2)$.
\end{proof}

Finally, lemmas~\ref{lemma:complexity2} and~\ref{lemma:complexity3}
yield Theorem~\ref{theo:second} (complexity result). The angular bijection 
ensures that the generating functions of $k$-rooted maps, $k$-rooted 
2-connected maps and $k$-rooted 3-connected maps are respectively equal 
to the series of $k$-rooted quadrangulations, $k$-rooted simple 
quadrangulations and $k$-rooted irreducible quadrangulations. Hence, 
lemmas~\ref{lemma:betarational},~\ref{lemma:etarational},~\ref{lemma:gammarational} and~\ref{lemma:2rootedirred} 
yield Theorem~\ref{theo:first} (algebraic structure).

\section{Conclusion}
We have proposed a new general and efficient method to enumerate
unrooted maps. In particular, we have improved significantly on the
complexity of counting oriented convex polyhedra (unrooted 3-connected
maps). 

Our method is flexible and can be adapted to enumerate other families
of unrooted maps. For example, a similar scheme can be used to count
unrooted loopless and then unrooted maps without loops and multiple
edges. This time, a first tree decomposition, ``by loops'', allows us
to enumerate $k$-rooted loopless maps from $k$-rooted
maps. Then the tree decomposition by multiple edges (this time on
$k$-rooted maps instead of $k$-rooted quadrangulations as in this
article) allows us to enumerate $k$-rooted maps without loop and multiple
edge from $k$-rooted loopless maps.

Another interesting problem is the enumeration of unrooted
3-connected maps on the sphere up to all homeomorphisms (including
orientation-reversing). Indeed according to Whitney's Theorem, 3-connected 
planar graphs have a unique toplogical embedding on the sphere, so that these
 unrooted 3-connected maps exactly correspond to unlabelled 3-connected 
planar graphs. In this case, another adaptation
of Burnside's
formula by Liskovets~\cite{Li96} is also available, giving an expression for the 
number of unrooted 
maps that involves the number of orientation-preserving $k$-rooted 
3-connected maps and also the number of orientation-reversing ones 
(for example 2-rooted 3-connected maps
where the 2-root induces a reflection). The tree-decomposition by separating 
4-cycles 
can be used to
obtain an equation linking 2-rooted 2-connected maps and 2-rooted
3-connected maps of type reflexion. Hence, the
method of tree decomposition is also here promising.

\appendix

\section{Enumeration with respect to the number of edges}
\label{sec:enum_edges}
\subsection{Unrooted maps}
\label{sec:enum_edges_1conn}
Let $\beta:=\beta(x)$ be the algebraic function defined by the equation:
$$
\beta(x)=x+3\beta(x)^2
$$
Let $f_n$ be the number of unrooted maps with $n$ edges. Then Burnside's formula for unconstrained maps is
$$
\sum_n 2nf_n x^n=F(x)+xF_{vf}(x^2)+x^2F_{f\!f}(x^2)+\sum_{k\geq 2}\phi(k)F_{vv}(x^k),
$$
where 
\begin{eqnarray*}
F(x)&=&{\frac {\beta\, \left( 2-9\,\beta \right) }{ \left( 1-3\,\beta
 \right) ^{2}}}\\
F_{vf}(x)&=&{\frac {2}{ \left( 1-6\,\beta \right)  \left( 1-3\,\beta \right) }}\\
F_{f\!f}(x)&=&{\frac {1}{ \left( 1-3\,\beta \right) ^{2} \left( 1-6\,\beta
 \right) }}\\
F_{vv}(x)&=&{\frac {6\beta}{1-6\,\beta}}.
\end{eqnarray*}
The first coefficients of the series of unrooted maps are $2\,x+4\,{x}^{2}+14\,{x}^{3}+57\,{x}^{4}+312\,{x}^{5}+2071\,{x}^{6}+15030\,{x}^{7}+117735\,{x}^{8}+967850\,{x}^{9}+8268816\,{x}^{10}+\ldots$

\subsection{Unrooted 2-connected maps}
\label{sec:enum_edges_2conn}
Let $\eta:=\eta(y)$ be the algebraic function defined by the equation:
$$
\eta(y)=y/(1-\eta(y))^2
$$
Let $g_n$ be the number of unrooted 2-connected maps with $n$ edges. Then Burnside's formula for 2-connected maps is
$$
\sum_n 2ng_n y^n=G(y)+yG_{vf}(y^2)+y^2G_{f\!f}(y^2)+\sum_{k\geq 2}\phi(k)G_{vv}(y^k),
$$
where 
\begin{eqnarray*}
G(y)&=&\eta\, \left( 2-3\,\eta \right) \\
G_{vf}(y)&=&\frac{2}{1-3\,\eta}\\
G_{f\!f}(y)&=&{\frac {1}{ \left( 1-3\,\eta \right)  \left( 1-\eta \right) }}\\
G_{vv}(y)&=&{\frac {2\eta}{1-3\,\eta}}.
\end{eqnarray*}
The first coefficients of the series of unrooted 2-connected maps are $2\,y+{y}^{2}+2\,{y}^{3}+3\,{y}^{4}+6\,{y}^{5}+16\,{y}^{6}+42\,{y}^{7}+151\,{y}^{8}+596\,{y}^{9}+2605\,{y}^{10}+\ldots$

\subsection{Unrooted 3-connected maps}
\label{sec:enum_edges_3conn}
Let $\gamma:=\gamma(z)$ be the algebraic function defined by the equation:
$$
\gamma(z)=z\left(1+\gamma(z)\right)^2
$$
Let $h_n$ be the number of unrooted 3-connected maps with $n$ edges. Then Burnside's formula for 3-connected maps is
$$
\sum_n 2nh_n z^n=H(z)+zH_{v\mathfrak{f}'}(z^2)+zH_{v\mathfrak{f}}(z^2)+z^2H_{\mathfrak{f}\mathfrak{f}'}(z^2)+z^2H_{\mathfrak{f}\mathfrak{f}}(z^2)+H_{vv}^{(2)}(z^2)+\sum_{k\geq 3}\phi(k)H_{vv}^{\geq 3}(z^k),
$$
where 
\begin{eqnarray*}
H(z)&=&-{\frac {{\gamma}^{6} \left( 2\,{\gamma}^{3}-1-4\,\gamma-3\,{\gamma}^{
2}+{\gamma}^{4} \right) }{ \left( 1+\gamma \right) ^{4} \left( 1+3\,
\gamma+{\gamma}^{2} \right) ^{2} \left( 2\,\gamma+1 \right) ^{3}}}\\
H_{v\mathfrak{f}}(z)&=&{\frac { 4\left( 1+\gamma \right)  \left( 8\,{\gamma}^{2}+13\,
\gamma+4 \right) {\gamma}^{4}}{ \left( 1-\gamma \right)  \left( 1+3\,
\gamma+{\gamma}^{2} \right) ^{2} \left( 2\,\gamma+1 \right) ^{3}}}
\\
H_{v\mathfrak{f}'}(z)&=&{\frac {2{\gamma}^{4}}{ \left( 1+3\,\gamma+{\gamma}^{2} \right) 
 \left( 2\,\gamma+1 \right) ^{2}}}
\\
H_{\mathfrak{f}\mathfrak{f}}(z)&=&{\frac {2{\gamma}^{2} \left( 1+5\,\gamma+10\,{\gamma}^{2}+9\,{
\gamma}^{3} \right)  \left( 1+\gamma \right) ^{2}}{ \left( 1-\gamma
 \right)  \left( 1+3\,\gamma+{\gamma}^{2} \right) ^{2} \left( 2\,
\gamma+1 \right) ^{3}}}
\\
H_{\mathfrak{f}\mathfrak{f}'}(z)&=&{\frac { \left( 1+3\,\gamma+3\,{\gamma}^{2} \right) {\gamma}^{2}}{
 \left( 1+3\,\gamma+{\gamma}^{2} \right)  \left( 2\,\gamma+1 \right) ^
{2}}}\\
H_{vv}^{(2)}(z)&=&{\frac { 2\left( 8\,{\gamma}^{5}+28\,{\gamma}^{4}+31\,{\gamma}^{3}+
21\,{\gamma}^{2}+10\,\gamma+2 \right) {\gamma}^{4}}{ \left( 1+3\,
\gamma+{\gamma}^{2} \right) ^{2} \left( 1+\gamma \right) ^{2} \left( 
1-\gamma \right)  \left( 2\,\gamma+1 \right) ^{3}}}\\
H_{vv}^{\geq 3}(z)&=&{\frac {2{\gamma}^{2} \left( 3\,\gamma+2 \right) }{ \left( 1-\gamma
 \right)  \left( 1+3\,\gamma+{\gamma}^{2} \right)  \left( 2\,\gamma+1
 \right) }}
\end{eqnarray*}
The first coefficients of the series of unrooted 3-connected maps are ${z}^{6}+{z}^{8}+2\,{z}^{9}+3\,{z}^{10}+4\,{z}^{11}+15\,{z}^{12}+32\,{z}^{13}+89\,{z}^{14}+266\,{z}^{15}+797\,{z}^{16}+2496\,{z}^{17}+\ldots$

\section{Enumeration with respect to the number of vertices and the number of faces}
\label{sec:enum_vertices}
\subsection{Unrooted maps}
\label{sec:enum_vertices_1conn}
Let $\beta_1:=\beta_1(\xb,\xw)$ and $\beta_2:=\beta_2(\xb,\xw)$ be the algebraic functions defined by the equation-system
$$
\left\{
\begin{array}{rcl}
\beta_1&=&\xb +\beta_1^2+2\beta_1\beta_2 \\
\beta_2&=&\xw +\beta_2^2+2\beta_1\beta_2.
\end{array}
\right.
$$
Let $f_{ij}$ be the number of unrooted maps with $i+1$ vertices and $j+1$ faces. Then Burnside's formula for unconstrained maps is
\begin{eqnarray*}
\sum_{i,j} 2(i+j)f_{ij} \xb^i\xw^j&=&F(\xb,\xw)+\xw F_{bf}(\xb^2,\xw^2)+\xb F_{wf}(\xb^2,\xw^2)+\xb\xw F_{f\!f}(\xb^2,\xw^2)\\
&&+\sum_{k\geq 2}\phi(k)\left(\frac{\xb}{\xw}F_{bb}(\xb^k,\xw^k)+\frac{\xw}{\xb}F_{ww}(\xb^k,\xw^k)+F_{bw}(\xb^k,\xw^k)\right),
\end{eqnarray*}
where 
\begin{eqnarray*}
F(\xb,\xw)&=&-{\frac {-\beta_2-\beta_1+5\,\beta_1\,\beta_2+2\,{\beta_1}^{2}+2\,{\beta_2}^{2}
}{ \left( -1+\beta_1+2\,\beta_2 \right)  \left( -1+\beta_2+2\,\beta_1
 \right) }}
\\
F_{bf}(\xb,\xw)&=&{\frac {-1+2\,\beta_2}{ \left( 4\,\beta_1\,\beta_2+1-4\,\beta_2-4\,\beta_1+4\,
{\beta_2}^{2}+4\,{\beta_1}^{2} \right)  \left( -1+\beta_2+2\,\beta_1
 \right) }}
\\
F_{wf}(\xb,\xw)&=&{\frac {-1+2\,\beta_1}{ \left( 4\,\beta_1\,\beta_2+1-4\,\beta_2-4\,\beta_1+4
\,{\beta_2}^{2}+4\,{\beta_1}^{2} \right)  \left( -1+\beta_1+2\,\beta_2
 \right) }}
\\
F_{f\!f}(\xb,\xw)&=&-{\frac {-1+\beta_1+\beta_2}{ \left( -1+\beta_1+2\,\beta_2 \right)  \left( -
1+\beta_2+2\,\beta_1 \right)  \left( 4\,\beta_1\,\beta_2+1-4\,\beta_2-4\,
\beta_1+4\,{\beta_2}^{2}+4\,{\beta_1}^{2} \right) }}
\\
F_{bb}(\xb,\xw)&=&{\frac {\beta_2}{4\,\beta_1\,\beta_2+1-4\,\beta_2-4\,\beta_1+4\,{\beta_2}^{2}+4
\,{\beta_1}^{2}}}
\\
F_{ww}(\xb,\xw)&=&{\frac {\beta_1}{4\,\beta_1\,\beta_2+1-4\,\beta_2-4\,\beta_1+4\,{\beta_2}^{2}+4
\,{\beta_1}^{2}}}
\\
F_{bw}(\xb,\xw)&=&-2\,{\frac {2\,{\beta_1}^{2}-\beta_1+2\,\beta_1\,\beta_2-\beta_2+2\,{\beta_2}^{
2}}{4\,\beta_1\,\beta_2+1-4\,\beta_2-4\,\beta_1+4\,{\beta_2}^{2}+4\,{\beta_1}^{
2}}}.
\\
\end{eqnarray*}
The first coefficients of the series of unrooted maps are $\left( \xw+\xb \right) + \left( {\xw}^{2}+2\,\xb\xw+{\xb}^{2} \right) + \left( 2\,{\xb}^{3}+5\,{\xb}^{2}\xw+5\,\xb{\xw}^{2}+2\,{\xw}^{3} \right) +
 \left( 14\,{\xb}^{3}\xw+14\,\xb{\xw}^{3}+23\,{\xb}^{2}{\xw}^{2}+3\,{\xw}^{4}+3\,{\xb}
^{4} \right) + \left( 108\,{\xb}^{2}{\xw}^{3}+6\,{\xb}^{5}+42\,{\xb}^{4
}\xw+6\,{\xw}^{5}+108\,{\xb}^{3}{\xw}^{2}+42\,\xb{\xw}^{4} \right) +\ldots$

\subsection{Unrooted 2-connected maps}
\label{sec:enum_vertices_2conn}
Let $\eta_1:=\eta_1(\yb,\yw)$ and $\eta_2:=\eta_2(\yb,\yw)$ be the algebraic functions defined by the equation-system
$$
\left\{
\begin{array}{rcl}
\eta_1&=&\yb/(1-\eta_2)^2\\
\eta_2&=&\yw/(1-\eta_1)^2.
\end{array}
\right.
$$
Let $g_{ij}$ be the number of unrooted 2-connected maps with $i+1$ vertices and $j+1$ faces. Then Burnside's formula for unrooted 2-connected maps is
\begin{eqnarray*}
\sum_{i,j} 2(i+j)g_{ij} \yb^i\yw^j&=&G(\yb,\yw)+\yw G_{bf}(\yb^2,\yw^2)+\yb G_{wf}(\yb^2,\yw^2)+\yb\yw G_{f\!f}(\yb^2,\yw^2)\\
&&+\sum_{k\geq 2}\phi(k)\left(\frac{\yb}{\yw}G_{bb}(\yb^k,\yw^k)+\frac{\yw}{\yb}G_{ww}(\yb^k,\yw^k)+G_{bw}(\yb^k,\yw^k)\right),
\end{eqnarray*}
where 
\begin{eqnarray*}
G(\yb,\yw)&=&-3\,\eta_1\,\eta_2+\eta_1+\eta_2
\\
G_{bf}(\yb,\yw)&=&-{\frac { \left( -1+\eta_2 \right)  \left( \eta_1+1 \right) }{ \left( -
1+\eta_1 \right)  \left( 3\,\eta_1\,\eta_2+\eta_2-1+\eta_1 \right) }}
\\
G_{wf}(\yb,\yw)&=&-{\frac { \left( \eta_2+1 \right)  \left( -1+\eta_1 \right) }{ \left( -
1+\eta_2 \right)  \left( 3\,\eta_1\,\eta_2+\eta_2-1+\eta_1 \right) }}
\\
G_{f\!f}(\yb,\yw)&=&{\frac {\eta_1\,\eta_2-1}{ \left( -1+\eta_2 \right)  \left( -1+\eta_1
 \right)  \left( 3\,\eta_1\,\eta_2+\eta_2-1+\eta_1 \right) }}
\\
G_{bb}(\yb,\yw)&=&{\frac {\eta_2\, \left( -1+\eta_1 \right) }{3\,\eta_1\,\eta_2+\eta_2-1+
\eta_1}}
\\
G_{ww}(\yb,\yw)&=&{\frac {\eta_1\, \left( -1+\eta_2 \right) }{3\,\eta_1\,\eta_2+\eta_2-1+
\eta_1}}
\\
G_{bw}(\yb,\yw)&=&-4\,{\frac {\eta_1\,\eta_2}{3\,\eta_1\,\eta_2+\eta_2-1+\eta_1}}
\\
\end{eqnarray*}
The first coefficients of the series of unrooted 2-connected maps are $\left( \yb+\yw \right) +\yb\yw+ \left( {\yb}^{2}\yw+\yb{\yw}^{2} \right)+ \left( {\yb}^{2}{\yw}^{2}+\yb{\yw}^{3}+{\yb}^{3}\yw \right)+\left( 2\,{\yb}^{3}{\yw}^{2}+2\,{\yb}^{2}{\yw}^{3}+{\yb}^{4}\yw+\yb{\yw}^{4} \right)+$\\
$\left( 3\,{\yb}^{2}{\yw}^{4}+3\,{\yb}^{4}{\yw}^{2}+\yb{\yw}^{5}+{\yb}^{5}\yw+8\,{\yb}^{3}{\yw}^{3} \right)+\ldots$

\subsection{Unrooted 3-connected maps}
\label{sec:enum_vertices_3conn}
Let $\gamma_1:=\gamma_1(\xb,\xw)$ and $\gamma_2:=\gamma_2(\xb,\xw)$ be the algebraic functions defined by the equation-system
$$
\left\{
\begin{array}{rcl}
\gamma_1&=&\zb(1+\gamma_2)^2\\
\gamma_2&=&\zw(1+\gamma_1)^2.
\end{array}
\right.
$$
Let $h_{ij}$ be the number of unrooted 3-connected maps with $i+1$ vertices and $j+1$ faces. Then Burnside's formula for unrooted 3-connected  maps is
\begin{eqnarray*}
\sum_{i,j} 2(i+j)h_{ij} \zb^i\zw^j&=&H(\zb,\zw)+\zw H_{b\mathfrak{f}'}(\zb^2,\zw^2)+\zw H_{b\mathfrak{f}}(\zb^2,\zw^2)+\zb H_{w\mathfrak{f}'}(\zb^2,\zw^2)\\
&&+\zb H_{w\mathfrak{f}}(\zb^2,\zw^2)+\zb\zw H_{\mathfrak{f}\mathfrak{f}'}(\zb^2,\zw^2)+\zb\zw H_{\mathfrak{f}\mathfrak{f}}(\zb^2,\zw^2)\\
&&+\frac{\zb}{\zw}H_{bb}^{(2)}(\zb^2,\zw^2)+\frac{\zw}{\zb}H_{ww}^{(2)}(\zb^2,\zw^2)+H_{bw}^{(2)}(\zb^2,\zw^2)\\
&&+\sum_{k\geq 3}\phi(k)\left(\frac{\zb}{\zw}H_{bb}^{\geq 3}(\zb^k,\zw^k)+\frac{\zw}{\zb}H_{ww}^{\geq 3}(\zb^k,\zw^k)+H_{bw}^{\geq 3}(\zb^k,\zw^k)\right),
\end{eqnarray*}
where 
{\small
\begin{eqnarray*}
H(\zb,\zw)&=&-{\frac {{\gamma_1}^{3}{\gamma_2}^{3} \left( {\gamma_1}^{2}{\gamma_2}^{2}-{
\gamma_1}^{2}+{\gamma_1}^{2}\gamma_2-\gamma_1\,\gamma_2-2\,\gamma_2-2\,\gamma_1-1-{
\gamma_2}^{2}+{\gamma_2}^{2}\gamma_1 \right) }{ \left( 1+\gamma_2 \right) ^{2}
 \left( 1+\gamma_1 \right) ^{2} \left( \gamma_1+1+2\,\gamma_2+{\gamma_2}^{2}
 \right)  \left( 1+2\,\gamma_1+{\gamma_1}^{2}+\gamma_2 \right)  \left( \gamma_1
+1+\gamma_2 \right) ^{3}}}
\\
H_{b\mathfrak{f}}(\zb,\zw)&=&-2\,{\gamma_1}^{2}{\gamma_2}^{2} \left( 4+16\,\gamma_2+13\,\gamma_1+24\,
{\gamma_2}^{2}+14\,{\gamma_1}^{2}+34\,\gamma_1\,\gamma_2+5\,{\gamma_1}^{3}+8\,{
\gamma_2}^{3}\gamma_1\right.\\
&&\left.+3\,{\gamma_1}^{3}\gamma_2+9\,{\gamma_1}^{2}{\gamma_2}^{2}+29\,{
\gamma_2}^{2}\gamma_1+16\,{\gamma_2}^{3}+4\,{\gamma_2}^{4}+21\,{\gamma_1}^{2}\gamma_2
 \right)/ \\ 
&& \left(\left( 3\,\gamma_1\,\gamma_2-\gamma_2-\gamma_1-1 \right)  \left( 
\gamma_1+1+2\,\gamma_2+{\gamma_2}^{2} \right) ^{2} \left( \gamma_1+1+\gamma_2
 \right) ^{3}\right)
\\
H_{w\mathfrak{f}}(\zb,\zw)&=&\mathrm{subs}(\{ \gamma_1=\gamma_2,\gamma_2=\gamma_1\},H_{b\mathfrak{f}})
\\
H_{b\mathfrak{f}'}(\zb,\zw)&=&{\frac {{\gamma_1}^{2}{\gamma_2}^{2}}{ \left( \gamma_1+1+2\,\gamma_2+{\gamma_2}^{2}
 \right)  \left( \gamma_1+1+\gamma_2 \right) ^{2}}}
\\
H_{w\mathfrak{f}'}(\zb,\zw)&=&\mathrm{subs}(\{ \gamma_1=\gamma_2,\gamma_2=\gamma_1\},H_{b\mathfrak{f}'})
%\end{eqnarray*}
%\begin{eqnarray*}
\\
H_{\mathfrak{f}\mathfrak{f}}(\zb,\zw)&=&
- \left( 2+112\,{\gamma_1}^{3}+112\,{\gamma_2}^{3}+56\,{\gamma_1}^{2}
+56\,{\gamma_2}^{2}+950\,{\gamma_1}^{2}{\gamma_2}^{2}+353\,{\gamma_1}^{2}\gamma_2+
114\,\gamma_1\,\gamma_2\right.\\
&&\left.+353\,{\gamma_2}^{2}\gamma_1+1708\,{\gamma_1}^{3}{\gamma_2}^{3
}+140\,{\gamma_2}^{4}+620\,{\gamma_2}^{3}\gamma_1+620\,{\gamma_1}^{3}\gamma_2+140\,
{\gamma_1}^{4}\right.\\
&&\left.+91\,{\gamma_1}^{3}{\gamma_2}^{6}+675\,{\gamma_1}^{4}\gamma_2+466\,{
\gamma_1}^{5}\gamma_2+683\,{\gamma_1}^{5}{\gamma_2}^{2}+10\,{\gamma_2}^{5}{\gamma_1}^
{5}+122\,{\gamma_2}^{5}{\gamma_1}^{4}\right.\\
&&\left.+122\,{\gamma_2}^{4}{\gamma_1}^{5}+445\,{
\gamma_2}^{3}{\gamma_1}^{5}+548\,{\gamma_1}^{4}{\gamma_2}^{4}+445\,{\gamma_1}^{3}{
\gamma_2}^{5}+1263\,{\gamma_2}^{2}{\gamma_1}^{4}+1164\,{\gamma_1}^{4}{\gamma_2}^{3}\right.\\
&&\left.+1415\,{\gamma_1}^{2}{\gamma_2}^{3}+1415\,{\gamma_1}^{3}{\gamma_2}^{2}+1263\,{
\gamma_1}^{2}{\gamma_2}^{4}+9\,{\gamma_1}^{3}{\gamma_2}^{7}+1164\,{\gamma_1}^{3}{
\gamma_2}^{4}+199\,{\gamma_1}^{6}\gamma_2\right.\\
&&\left.+216\,{\gamma_1}^{6}{\gamma_2}^{2}+91\,{
\gamma_1}^{6}{\gamma_2}^{3}+10\,{\gamma_1}^{6}{\gamma_2}^{4}+37\,{\gamma_1}^{7}{
\gamma_2}^{2}+9\,{\gamma_1}^{7}{\gamma_2}^{3}+48\,{\gamma_1}^{7}\gamma_2+5\,{\gamma_1
}^{8}\gamma_2\right.\\
&&\left.+3\,{\gamma_1}^{8}{\gamma_2}^{2}+675\,\gamma_1\,{\gamma_2}^{4}+466\,
\gamma_1\,{\gamma_2}^{5}+199\,\gamma_1\,{\gamma_2}^{6}+48\,\gamma_1\,{\gamma_2}^{7}+5
\,\gamma_1\,{\gamma_2}^{8}\right.\\
&&\left.+683\,{\gamma_1}^{2}{\gamma_2}^{5}+216\,{\gamma_1}^{2}{
\gamma_2}^{6}+37\,{\gamma_1}^{2}{\gamma_2}^{7}+3\,{\gamma_1}^{2}{\gamma_2}^{8}+10\,
{\gamma_1}^{4}{\gamma_2}^{6}+2\,{\gamma_1}^{8}+56\,{\gamma_1}^{6}\right.\\
&&\left.+16\,{\gamma_1}^
{7}+112\,{\gamma_1}^{5}+56\,{\gamma_2}^{6}+2\,{\gamma_2}^{8}+16\,{\gamma_2}^{7}+
112\,{\gamma_2}^{5}+16\,\gamma_1+16\,\gamma_2 \right) \gamma_1\,\gamma_2\ / \\
&& \left(\left( 3
\,\gamma_1\,\gamma_2-\gamma_2-\gamma_1-1 \right)  \left( \gamma_1+1+2\,\gamma_2+{
\gamma_2}^{2} \right) ^{2} \left( 1+2\,\gamma_1+{\gamma_1}^{2}+\gamma_2 \right) 
^{2} \left( \gamma_1+1+\gamma_2 \right) ^{3}\right)
%\\
\end{eqnarray*}
\begin{eqnarray*}
H_{\mathfrak{f}\mathfrak{f}'}(\zb,\zw)&=&\left( 1+3\,\gamma_2+3\,\gamma_1+3\,{\gamma_2}^{2}+3\,{\gamma_1}^{2}+7
\,\gamma_1\,\gamma_2+{\gamma_1}^{3}+{\gamma_2}^{3}\gamma_1+{\gamma_1}^{3}\gamma_2+{
\gamma_1}^{2}{\gamma_2}^{2}\right.\\
&&\left.+5\,{\gamma_2}^{2}\gamma_1+{\gamma_2}^{3}+5\,{\gamma_1}^{2
}\gamma_2 \right) \gamma_1\,\gamma_2\ / \\
&&\left(\left( 1+2\,\gamma_1+{\gamma_1}^{2}+\gamma_2
 \right)  \left( \gamma_1+1+2\,\gamma_2+{\gamma_2}^{2} \right)  \left( \gamma_1+
1+\gamma_2 \right) ^{2}\right)
\\
H_{bb}^{(2)}(\zb,\zw)&=&-{\gamma_1}^{2}{\gamma_2}^{4} \left( 3+15\,\gamma_2+20\,\gamma_1+30\,
\gamma_1\,{\gamma_2}^{4}+3\,\gamma_1\,{\gamma_2}^{5}+11\,{\gamma_1}^{4}\gamma_2+12\,{\gamma_1}^{4}+30\,{\gamma_2}^{2}\right.\\
&&\left.+3\,{\gamma_2}^{2}{\gamma_1}^{4}+43\,{\gamma_1}^{2
}+9\,{\gamma_1}^{3}{\gamma_2}^{3}+81\,\gamma_1\,\gamma_2+55\,{\gamma_1}^{2}{\gamma_2}
^{3}+38\,{\gamma_1}^{3}+43\,{\gamma_1}^{3}{\gamma_2}^{2}\right.\\
&&\left.+92\,{\gamma_2}^{3}
\gamma_1+72\,{\gamma_1}^{3}\gamma_2+132\,{\gamma_1}^{2}{\gamma_2}^{2}+126\,{\gamma_2}
^{2}\gamma_1+3\,{\gamma_2}^{5}+30\,{\gamma_2}^{3}+15\,{\gamma_2}^{4}\right.\\
&&\left.+127\,{\gamma_1
}^{2}\gamma_2+7\,{\gamma_1}^{2}{\gamma_2}^{4} \right)\ / \\
&& \left(\left( 3\,\gamma_1\,
\gamma_2-\gamma_2-\gamma_1-1 \right)  \left( 1+2\,\gamma_1+{\gamma_1}^{2}+\gamma_2
 \right)  \left( 1+\gamma_1 \right) ^{2} \left( \gamma_1+1+2\,\gamma_2+{\gamma_2
}^{2} \right) ^{2}\right.\\
&& \left.\left( \gamma_1+1+\gamma_2 \right) ^{3}\right)
\\
H_{ww}^{(2)}(\zb,\zw)&=&\mathrm{subs}(\{ \gamma_1=\gamma_2,\gamma_2=\gamma_1\},H_{bb}^{(2)})
\\
H_{bw}^{(2)}(\zb,\zw)&=&-4\,{\frac {{\gamma_1}^{2}{\gamma_2}^{2}}{ \left( 3\,\gamma_1\,\gamma_2-\gamma_2-
\gamma_1-1 \right)  \left( \gamma_1+1+\gamma_2 \right) ^{3}}}
\\
H_{bb}^{\geq 3}(\zb,\zw)&=&-{\frac {{\gamma_2}^{2}\gamma_1\, \left( 4\,\gamma_1+3+3\,\gamma_2 \right) }{
 \left( 3\,\gamma_1\,\gamma_2-\gamma_2-\gamma_1-1 \right)  \left( 1+2\,\gamma_1+{
\gamma_1}^{2}+\gamma_2 \right)  \left( \gamma_1+1+\gamma_2 \right) }}
\\
H_{ww}^{\geq 3}(\zb,\zw)&=&\mathrm{subs}(\{ \gamma_1=\gamma_2,\gamma_2=\gamma_1\},H_{bb}^{\geq 3})
\\
H_{bw}^{\geq 3}(\zb,\zw)&=&-4\,{\frac {\gamma_1\,\gamma_2}{ \left( \gamma_1+1+\gamma_2 \right)  \left( 3\,
\gamma_1\,\gamma_2-\gamma_2-\gamma_1-1 \right) }}
\\
\end{eqnarray*}}
The first coefficients of the series of unrooted 3-connected maps are ${\zb}^{3}{\zw}^{3}+{\zb}^{4}{\zw}^{4}+ \left( {\zb}^{5}{\zw}^{4}+{\zb}^{4}{\zw}^{5} \right)+3\,{\zb}^{5}{\zw}^{5}+ \left( 2\,{\zb}^
{6}{\zw}^{5}+2\,{\zb}^{5}{\zw}^{6} \right) + \left( 2\,{\zb}^{7}{\zw}^{5
}+2\,{\zb}^{5}{\zw}^{7}+11\,{\zb}^{6}{\zw}^{6} \right)+ \left( 16\,{\zb
}^{7}{\zw}^{6}+16\,{\zb}^{6}{\zw}^{7} \right)+ \left( 10\,{\zb}^{8}{\zw
}^{6}+69\,{\zb}^{7}{\zw}^{7}+10\,{\zb}^{6}{\zw}^{8} \right) +\ldots$

\end{document}

%% file: Figures/head.tex
\usepackage{calc} % pour \overpar (parenthese au dessus d'une expression pour \ronde{})
\usepackage{amsmath}
\usepackage{amsfonts}
\usepackage{array}
\usepackage[latin1]{inputenc}
\usepackage[mathscr]{eucal} % définit \mathscr{} (jolie fonte)

\renewcommand{\leq}{\ensuremath{\leqslant}}
\renewcommand{\geq}{\ensuremath{\geqslant}}

\renewcommand{\epsilon}{\ensuremath{\varepsilon}}

\newlength{\restsubwidth}
\newlength{\restsubheight}
\newlength{\restsubmoreheight}
\newcommand{\rest}[2]{%
    \settowidth{\restsubwidth}{\ensuremath{#2}}
    \settoheight{\restsubheight}{\ensuremath{{}_{#2}}}
    \ensuremath{{#1\hskip 0.5 pt}_{\vrule\kern2pt\parbox[b][%
    \the\restsubheight +
        \the\restsubmoreheight][b]{\the\restsubwidth}{%
            \ensuremath{{}_{#2}}}}}
    }

\renewcommand{\vec}[1]{\vect{#1}}

\mathcode`A="7041 \mathcode`B="7042 \mathcode`C="7043 \mathcode`D="7044
\mathcode`E="7045 \mathcode`F="7046 \mathcode`G="7047 \mathcode`H="7048
\mathcode`I="7049 \mathcode`J="704A \mathcode`K="704B \mathcode`L="704C
\mathcode`M="704D \mathcode`N="704E \mathcode`O="704F \mathcode`P="7050
\mathcode`Q="7051 \mathcode`R="7052 \mathcode`S="7053 \mathcode`T="7054
\mathcode`U="7055 \mathcode`V="7056 \mathcode`W="7057 \mathcode`X="7058
\mathcode`Y="7059 \mathcode`Z="705A

\newlength{\moreinterligne}
\newlength{\letterarrowvskip}
\newlength{\vectlength}
\newlength{\vectheight}
\newlength{\boxlength}
\newcommand{\vect}[1]{%
        \settowidth{\vectlength}{\ensuremath{#1}}%
        \settoheight{\vectheight}{\ensuremath{#1}}%
        \settowidth{\boxlength}{%
        \ensuremath{%
        \overrightarrow{%
        \parbox[b][\the\vectheight + %
                \letterarrowvskip][b]{\the\vectlength}%
        {\ensuremath{#1}}}}%
        }
        \parbox[b][\the\vectheight
                + \moreinterligne][b]{\the\boxlength}{%
        \ensuremath{%
    \overset{\hbox to \the\boxlength{\rightarrowfill}}
    {
        \parbox[b][\the\vectheight + %
                \letterarrowvskip][b]{\the\vectlength}%
        {\ensuremath{#1}}}}%
        }
        }

\newlength{\smoreinterligne}
\newlength{\sletterarrowvskip}
\newlength{\svectlength}
\newlength{\svectheight}
\newlength{\sboxlength}
\newcommand{\svect}[1]{%
        \settowidth{\svectlength}{\mbox{\scriptsize{\ensuremath{#1}}}}%
        \settoheight{\svectheight}{\mbox{\scriptsize{\ensuremath{#1}}}}%
        \settowidth{\sboxlength}{%
        \ensuremath{%
        \overrightarrow{%
        \parbox[b][\the\svectheight + %
                \sletterarrowvskip][b]{\the\svectlength}%
        {\mbox{\scriptsize{\ensuremath{#1}}}}}}%
        }
        \parbox[b][\the\svectheight
                + \smoreinterligne][b]{\the\sboxlength}{%
        \ensuremath{%
        \overrightarrow{%
        \parbox[b][\the\svectheight + %
                \sletterarrowvskip][b]{\the\svectlength}%
        {\mbox{\scriptsize{\ensuremath{#1}}}}}}%
        }
        }

\newlength{\moreinterligneo}
\newlength{\letterarrowvskipo}
\newlength{\vectlengtho}
\newlength{\vectheighto}
\newlength{\boxlengtho}
\newcommand{\vecttitre}[1]{%
        \settowidth{\vectlengtho}{\ensuremath{#1}}%
        \settoheight{\vectheighto}{\ensuremath{#1}}%
        \settowidth{\boxlengtho}{%
        \ensuremath{%
        \overrightarrow{%
        \parbox[b][\the\vectheighto + %
                \letterarrowvskipo][b]{\the\vectlengtho}%
        {\ensuremath{#1}}}}%
        }
        \parbox[b][\the\vectheighto
                + \moreinterligneo][b]{\the\boxlengtho}{%
        \ensuremath{%
    \overset{\hbox to \the\boxlengtho{\rightarrowfill}}
    {
        \parbox[b][\the\vectheighto + %
                \letterarrowvskipo][b]{\the\vectlengtho}%
        {\ensuremath{#1}}}}%
        }
        }

\newlength{\remlength}

    \newlength{\TrigArgUnH}
    \newlength{\TrigArgUnV}
    \newlength{\TrigArgDeuxH}
    \newlength{\TrigArgDeuxV}
    \newlength{\tempH}
    \newlength{\tempV}
    \newlength{\TrigArgTroisH}
    \newlength{\TrigCentreH}
    \newlength{\TrigCentreV}
    \newlength{\decalage}
\newenvironment{ltab}[1][2]{%
        }{%
        }

\newlength{\espace} 
\newlength{\letterparvskip} 
\newlength{\parlength}
\newlength{\parheight}
\newlength{\arglength}
\newlength{\complength}
\newlength{\lmouslength}
\newlength{\rmouslength}
\newcommand{\overpar}[1]{%
        \settowidth{\parlength}{\ensuremath{#1}}%
        \settoheight{\parheight}{\ensuremath{#1}}%
        \settowidth{\arglength}{%
        \ensuremath{%
        \overrightarrow{%
        \parbox[b][\the\parheight + %
                \letterparvskip][b]{\the\parlength}%
        {\ensuremath{#1}}}}%
        }
        \parbox[b][\the\parheight
                + \espace][b]{\the\arglength}{%
        \setlength{\complength}{\the\arglength - \the\lmouslength%
                - \the\rmouslength}
        \ensuremath{%
        \overset{\lmoustache\rule{\the\complength}{1 pt}\rmoustache}
        {
        \parbox[b][\the\parheight + %
                \letterparvskip][b]{\the\parlength}%
        {\ensuremath{#1}}}%
        }
        }
        }

\newlength{\alignelength}

\def\leaderfill{\leaders\hbox to 1ex{\hss.\hss}\hfill}

%% file: Figures/declh.pstex_t
\begin{picture}(0,0)%
\includegraphics{Figures/declh.pstex}%
\end{picture}%
\setlength{\unitlength}{2289sp}%
\begingroup\makeatletter\ifx\SetFigFont\undefined%
\gdef\SetFigFont#1#2#3#4#5{%
  \reset@font\fontsize{#1}{#2pt}%
  \fontfamily{#3}\fontseries{#4}\fontshape{#5}%
  \selectfont}%
\fi\endgroup%
\begin{picture}(9020,2512)(621,-2240)
\put(1651,-436){\makebox(0,0)[lb]{\smash{{\SetFigFont{10}{12.0}{\rmdefault}{\mddefault}{\updefault}{\color[rgb]{0,0,0}$\Gfftt$}%
}}}}
\put(1651,-1861){\makebox(0,0)[lb]{\smash{{\SetFigFont{10}{12.0}{\rmdefault}{\mddefault}{\updefault}{\color[rgb]{0,0,0}$\Gfvtt$}%
}}}}
\end{picture}%

%% file: Figures/dech.pstex_t
\begin{picture}(0,0)%
\includegraphics{Figures/dech.pstex}%
\end{picture}%
\setlength{\unitlength}{2644sp}%
\begingroup\makeatletter\ifx\SetFigFont\undefined%
\gdef\SetFigFont#1#2#3#4#5{%
  \reset@font\fontsize{#1}{#2pt}%
  \fontfamily{#3}\fontseries{#4}\fontshape{#5}%
  \selectfont}%
\fi\endgroup%
\begin{picture}(9140,1846)(1226,-1686)
\put(2006,-551){\makebox(0,0)[lb]{\smash{{\SetFigFont{12}{14.4}{\familydefault}{\mddefault}{\updefault}{\color[rgb]{0,0,0}$\mathcal{K}_{b}$}%
}}}}
\put(7911,-569){\makebox(0,0)[lb]{\smash{{\SetFigFont{12}{14.4}{\rmdefault}{\mddefault}{\updefault}{\color[rgb]{0,0,0}$\Gfbtt$}%
}}}}
\put(6075,-270){\makebox(0,0)[lb]{\smash{{\SetFigFont{12}{14.4}{\rmdefault}{\mddefault}{\updefault}{\color[rgb]{0,0,0}\rotatebox{270}{$^t\mathcal{K}_{w}$}}%
}}}}
\put(9491,-564){\makebox(0,0)[lb]{\smash{{\SetFigFont{12}{14.4}{\rmdefault}{\mddefault}{\updefault}{\color[rgb]{0,0,0}$\Qfbtt$}%
}}}}
\end{picture}%